\documentclass[12pt]{article}

\usepackage{amsthm,amsmath,amssymb,amsfonts}
\usepackage{epsfig,graphicx,float,subfig}

\usepackage{array}
\usepackage{listings}
\usepackage{wrapfig}
\usepackage{enumerate}

\theoremstyle{plain}
\newtheorem{theorem}{Theorem}
\newtheorem{corollary}{Corollary} 

\usepackage{graphicx, epsfig, amssymb}

\newcommand{\Prob} {{\bf P}}

\newcommand{\E}{{\bf E}}

\begin{document}

\title{Large-Scale Multi-Stream Quickest Change Detection via Shrinkage Post-Change Estimation}

\author{Yuan~Wang~and~Yajun~Mei
\thanks{The
material in this article was presented in part at the IEEE
International Symposium on Information Theory, Honolulu, HI, USA,
2014.
}%
\thanks{The authors are with The H. Milton Stewart School of Industrial and Systems
Engineering at the Georgia Institute of Technology, Atlanta,
Georgia, USA ~(email: sxtywangyuan@gmail.com, ymei@isye.gatech.edu).}
\thanks{Copyright (c) 2014 IEEE} }

\date{July 10, 2015}
\maketitle

\begin{abstract}
The quickest change detection problem is considered in the context of monitoring  large-scale  independent normal distributed data streams with possible changes in some of the means. It is assumed that  for each individual local data stream, either there are no local changes, or there is a ``big" local change that is larger than a pre-specified lower bound. Two different kinds of  scenarios are studied: one is the sparse post-change case when the unknown number of affected data streams is much smaller than the total number of data streams, and the other is when all local data streams are affected simultaneously although not necessarily identically. We propose a systematic approach to develop efficient global monitoring schemes for quickest change detection by combining hard thresholding with linear shrinkage estimators to estimating all post-change parameters  simultaneously. Our theoretical analysis  demonstrates that the shrinkage estimation can balance the tradeoff between the first-order and second-order terms of the asymptotic expression on the detection delays, and our numerical simulation studies illustrate the usefulness of shrinkage estimation and the challenge of Monte Carlo simulation of the average run length to false alarm in the context of online monitoring large-scale data streams.
\end{abstract}


\section{Introduction}

The  problem of online monitoring large-scale data streams
has many important applications from biosurveillance and quality control to
finance and security in modern information age when the rapid development
of sensing technology allows one to generate large-scale real-time streaming data.
In many scenarios, one is often interested in the early detection
of a ``trigger" event when ``sensors" are deployed to monitor
the changing environments over time and space, see
Lawson and Kleinman \cite{lawson2005spatial}. From the theoretical or methodological viewpoint,
this is a quickest change detection or sequential change-point detection problem,
where the case of monitoring $K=1$ data stream has been extensively studied in
the past several decades, see the books by Basseville and Nikiforov \cite{basseville1993detection} and Poor and Hadjiliadis \cite{poor2009quickest} for the review.
Also see Page \cite{page1954continuous}, Shiryaev \cite{shiryaev1963optimum}, Lorden \cite{lorden1971procedures}, Pollak \cite{pollak1985optimal}, Moustakides \cite{moustakides1986optimal}, Lai \cite{lai1995sequential} for some early classical contributions. In addition, the case of online monitoring a not so large number  $K$ (e.g., tens) of data streams
has also been studied in the literature, see Lorden and Pollak \cite{lorden2005nonanticipating}, Tartakovsky et al. \cite{tartakovsky2006detection}, Zamba and Hawkins \cite{zamba2006multivariate}, Veeravalli and Bangerjee \cite{veeravalli2012quickest}.

Unfortunately research on the problem of online monitoring a large  number $K$ (e.g., hundreds or more) of data streams is rather limited,  see Siegmund \cite{siegmund2013change} and the discussions therein. 
While many classical quickest change detection methods are based on likelihood ratio statistics, and can be extended from one-dimensional to $K$-dimensional, their performances are rather poor when monitoring a large number $K$ of data streams, despite holding the so-called first-order asymptotic optimality properties for any fixed dimensional $K$ in the sense of asymptotically minimizing the detection delay for each and every possible post-change hypothesis as the average run length (ARL)  to false alarm constraint goes to $\infty.$ The main reason is that these classical quickest change detection methods often over-emphasize the first-order performance for {\it each and every possible} post-change hypothesis in the $K$-dimensional space, and thus the price they paid is on the second-order terms of the detection delays which are often linearly increasing as a function of $K.$ This is not an issue when the number $K$ of data streams is small, but it has a severe effect when $K$ is large (e.g., hundreds): under a reasonable practical setting, the second-order term of the detection delay will likely be comparable to the first-order term, which implies that the nice first-order asymptotic optimality properties have little practical meaning for large $K$! This led Mei \cite{mei2010efficient} to raise an open problem whether one can develop new methods that can reduce the coefficient in the second-order term of the detection delay from  $K$  to a smaller number to yield quicker detection.

The primary objective of this paper is to tackle this open problem, and propose a systematic approach to develop efficient methodologies for online monitoring a large number $K$ of independent data streams. Our proposed methods do not aim for each and every possible post-change hypothesis in the $K$-dimensional space, and the main assumption we make is that for each individual local data stream, either there are no local changes, or there is a local change that is larger than some pre-specified lower bounds. The key novelty of our proposed methodologies  is to apply shrinkage estimators to incorporate such prior knowledge of the post-change hypothesis to develop efficient quickest change detection methodologies.  To illustrate our main ideas, we will focus on the problem of monitoring $K$ independent normal data streams with possible changes in the means of some data streams, and two different scenarios will be investigated: one is the sparse post-change case when the unknown number of affected data streams is much smaller than the total number of data streams, and the other is when all local data streams are affected simultaneously although not necessarily identically, i.e., different local data streams may have different {\it unknown} post-change mean  parameters.  It is useful to think that for a given total information for changing event, the former scenario corresponds to the case of a few ``large" local changes, whereas the latter scenario corresponds to the case of ``relatively small" local changes in all data streams. Given the same total information of changing event, the classical quickest change detection methods will have similar (first-order) performance under these two scenarios, although intuitively one may feel that these two scenarios should be different. Our proposed methods combine the hard thresholding estimators with the linear shrinkage estimators to simultaneously estimate unknown post-change mean parameters, and will indeed show that these two scenarios should be treated differently.
In the process of investigating the properties of the proposed methods, we also demonstrate the challenge of Monte Carlo simulation of the average run length to false alarm for large dimensional $K$ due to the curse of dimensionality, which seems to be overlooked in the quickest change detection literature. 

Note that the usefulness of shrinkage or thresholding in high-dimensional data is well-known in the modern off-line statistical research since the pioneering work of James and Stein \cite{james1961estimation}, also see Cand\'{e}s \cite{candes2006modern} and references therein. However, the application of shrinkage or thresholding to quickest change detection is rather limited.  Unlike other off-line works that deal with high-dimensional statistics, the asymptotic analysis in this paper fixes the dimension $K$ (or the number of data streams) as the ARL to false alarm is taken to infinity. Our aim is on the development of asymptotic results that are useful for the practical setting, and thus our focuses are on the effects of the dimension $K$ on the second-order term of the detection delays, and on how shrinkage or thresholding can lessen such effects.  As far as we know, it remains an open problem in quickest change detection  when the dimension $K$ is taken to infinity.

In the present paper, we will demonstrate how to combine shrinkage estimators with the classical Shiryaev-Roberts procedure to yield an efficient global monitoring scheme.
Note that the Shiryaev-Roberts procedure is chosen as a demonstration here, since it allows us to simplify our mathematical arguments by borrowing the results in Lorden and Pollak \cite{lorden2005nonanticipating} that develops the Shiryaev-Roberts-Robbins-Siegmund (SRRS) scheme based on the method of moments (MOM) estimators or the maximum likelihood estimators (MLE) of unknown post-change parameters. Besides the different estimators of unknown post-change parameters, another main difference between our research and Lorden and Pollak \cite{lorden2005nonanticipating} is that we explicitly investigate the effect of the number of data streams on the detection delay performance of the schemes.
We want to emphasize that our use of shrinkage estimators can easily be combined to other popular quickest change methods such as the CUSUM procedure proposed by Page \cite{page1954continuous} from the methodology or algorithm point of view, although the corresponding theoretical asymptotic analysis seems to be nontrivial. Hopefully our useful of shrinkage estimation opens new directions to develop more efficient methodologies for online monitoring of large-scale or high-dimensional data streams.

From the information theory viewpoint, the asymptotic performance of our proposed shrinkage-based schemes is characterized by the new information number  defined in (\ref{eqn02sec5-info_gen})  below. In a simple setting for normal distributions when the $\omega_{k}$'s are the smallest meaningful bounds on the post-change mean parameters $\mu_{k}$'s, the new information number has the form of $\frac12 \sum_{k: |\mu_{k}| > \omega_{k}} (\mu_{k})^2,$ whereas the classical Kullback-Leibler divergence is $\frac12 \sum_{k=1}^{K} (\mu_{k})^2.$  Thus our proposed new information number can be thought of as the shrinkage approximation of the classical Kullback-Leibler divergence  between pre-change and post-change distributions. In the context of monitoring large-scale data streams, we feel that our proposed new information number in  (\ref{eqn02sec5-info_gen}) provides more meaningful bounds than the  classical Kullback-Leibler divergence, since it takes into account of  the second-order term of the detection delay performance and the spatial uncertainty associated with which local data streams are affected. 

We should acknowledge that Xie and Siegmund \cite{xie2013sequential} studies a similar problem by taking a semi-Bayesian approach under the assumption that the fraction of affected data streams is known. Here we did not make such an assumption, and our formulation assumes that the lower bound of the post-change parameters are given, i.e., we are only interested in detecting certain large local changes for individual local data streams. In addition, Tartakovsky et al. \cite{tartakovsky2006detection} and Mei \cite{mei2010efficient} consider the special case when all post-change parameters for affected data streams were identical or completely specified. Here our underlying assumption is that the post-change parameters are unknown and not necessarily identical. In addition, the problem of monitoring $K>1$ data streams is also studied in the offline setting when the full information is available during decision-making, e.g. Zhang, Siegmund, Ji and Li \cite{zhang2010detecting}, and Cho and Fryzlewicz \cite{cho2012multiple}. Our setting here is online where we observe the data sequentially over time, and we cannot use future observations to make current   decision.


The remainder of this paper is organized as follows. In Section II, we state the mathematical formulation of monitoring $K > 1$ data streams and review shrinkage estimators in offline point estimation that will be used later. In Section III, we propose our shrinkage-based monitoring scheme for the problem of online monitoring of independent normal data streams with possible changes in some of the means. Section IV develops asymptotic properties of our proposed monitoring schemes. In Section V, we report numerical simulation results to illustrate the usefulness of our proposed shrinkage-based schemes  and  the challenge of Monte Carlo simulation of the average run length to false alarm in the context of online monitoring large-scale data streams.    Section VI contains some concluding remarks. 
The proof of Theorem \ref{thm: detection change} is included in the Appendix.

\section{Problem formulation and Background}

\subsection{Problem formulation}

Assume we are monitoring $K$ independent normal data streams in a system. Denote by $X_{k,n}$ the observation of the $k$-th data stream at time $n$ for $k = 1,...,K$ and $n=1,2, ... .$  The $X_{k,n}$'s are assumed to be independent not only over time within each data stream, but also among different data streams. Initially, all $X_{k,n}$'s are independent and identically distributed (iid) $N(\mu_0,1)$ random variables.  At some unknown time $\nu \in \{1,2,3,\ldots\},$ an event may occur to the system, and affect some data streams in the sense that
the distribution of the $X_{k,n}$'s may change to $N(\mu_k, 1)$ for $n = \nu, \nu+1, ... ,$ if the $k$-th data stream  is affected for $k=1, \ldots, K.$ To simplify our notation, here the post-change mean $\mu_{k} = \mu_0$ implies that the corresponding data stream is not affected, whereas $\mu_{k} \ne \mu_0$ corresponds to an affected data stream.  Following the literature of quickest change detection, we assume that the pre-change mean $\mu_0$ is completely specified, and without loss of generality, we assume $\mu_0 = 0,$ as otherwise we can monitor $X_{k,n} - \mu_0$ instead of $X_{k,n}$'s themselves. Thus $\mu_0$ and $0$ are interchangeable below for normal distributions.


In this article, we tackle the case when the post-change means $\mu_{k}$'s are only partially specified, e.g., we do not know which data streams are affected and do not know the exact values of the post-change means $\mu_{k}$'s for affected data streams. In practical situations when monitoring large-scale data streams,  one is often interested in only detecting ``big" local changes in individual data streams. This motivates us to assume that the post-change hypothesis set for the  post-change mean vector $\mu = (\mu_1, \ldots, \mu_K)^{T}$ is given by
\begin{eqnarray} \label{par:par}
\Omega =\{\mu \ne {\bf 0} \in \mathcal{R}^{K}: \sum_{k=1}^{K}  |\mu_{k}| 1\{|\mu_k| \le \omega_{k}\} = 0  \},
\end{eqnarray}
where the lower bounds $\omega_{k}$'s are pre-specified positive constants that are the smallest difference meaningful for detection. The post-change hypothesis set $\Omega$ in  (\ref{par:par}) implies that for any local data stream, either there are no local changes (i.e., $\mu_{k} = 0$), or there is a big local change (i.e., $|\mu_{k}| > \omega_{k}$). In addition, $\mu \ne {\bf 0}$ implies that at least  one $\mu_{k} \ne 0,$ i.e., at least one data stream should be affected under the post-change hypothesis.   Also note that the post-change hypothesis set $\Omega$ in  (\ref{par:par}) assumes  the true post-change mean $\mu_{k} \ne \pm \omega_{k}$ for any $k.$ This is a technical assumption to simplify our theoretical analysis, since otherwise  careful arguments are needed to take care of those data streams with $|\mu_{k}| = \omega_{k} > 0$ which could be thought of as affected data streams only with probability $1/2.$ For any given post-change mean vector $\mu,$  it is natural to define the number of affected data streams as  $r = \sum_{k=1}^{K} 1\{ \mu_k \ne 0 \}$ where $1\{A\}$ is the indicator function of event $A.$ Clearly, when $\mu  \in \Omega$ in (\ref{par:par}), this becomes
\begin{eqnarray} \label{par:par3}
r =  \sum_{k=1}^{K} 1\{ |\mu_k| > \omega_{k} \},
\end{eqnarray}
which will play an important role on the detection delay performance of quickest change detection schemes in our context. Note that the main scheme in Xie and Siegmund \cite{xie2013sequential} assumes that  the number of affected data streams, or the $r$ value in (\ref{par:par3}), is known and the lower bound $\omega_{k} = 0$  for all $k.$ In this article, we assume that the lower bounds $\omega_{k}$'s in (\ref{par:par})  are known {\bf positive} constants for all $k=1, \ldots, K.$ Two scenarios will be studied: one is the sparse post-change hypothesis case when the value $r$ in (\ref{par:par3}) is an unknown constant that is much smaller than $K,$ and the other is when $r= K,$ i.e., when all data streams are affected simultaneously.

To provide a more rigorous mathematical formulation, denote by $\Prob_{\mu,\nu}$ and $\E_{\mu,\nu}$ the probability measure and expectation of $\{(X_{k,1}$, $X_{k,2},$$...)\}_{k=1}^{p}$ when the change occurs at time $\nu$ and the true post-change mean vector $\mu = (\mu_1, \ldots, \mu_p)^{T}$. Denote by $\Prob_{\infty}$ and $\E_{\infty}$ the same when no change occurs, i.e., the change-time $\nu = \infty.$ Loosely speaking, we want to develop an online global monitoring scheme that can raise a true alarm as soon as possible when the event occurs while controlling the global false alarm rate. Mathematically, an online global monitoring scheme is defined as a stopping time $T,$ which is an integer-valued random variable. The event $\{T=n\}$ represents that we will raise an alarm at time $n$ at the global level and declare that a change occurs somewhere in the first $n$ time steps. Note that the decision $\{T=n\}$ is only based on the observations $X_{k,i}$'s up to time $n$.

The standard minimax formulation of quickest change detection problem can then be formally stated as follows: Find  a stopping time $T$ that asymptotically minimizes the ``worst-case" detection delay proposed in Lorden \cite{lorden1971procedures}
\begin{eqnarray*}
D_{\mu}(T) =  \sup_{1 \leq \nu < \infty} \mbox{ess sup}\ \E_{\mu, \nu} (T - \nu + 1 | T \geq \nu, \mathcal{F}_{\nu-1})
\end{eqnarray*}
for all possible post-change mean vectors $\mu \in \Omega$ in (\ref{par:par})  subject to the constraint on the average run length (ARL) to false alarm
\begin{eqnarray}
    \E_{\infty} (T) \geq A.                         \label{eq: change detect cond}
\end{eqnarray}
Here $\mathcal{F}_{\nu-1}$ denotes all information up to time $\nu-1,$ and the constraint $A > 0$ in  (\ref{eq: change detect cond}) is pre-specified.

\subsection{Review of Shrinkage Estimation}

Let us now review some well-known fact regarding offline shrinkage estimation procedures, which will be used in our proposed methodologies for online monitoring $K > 1$ data streams in the next section. Suppose that there are $K \ge 3$ independent normal random variables, say, $\{Y_1, \ldots, Y_{K}\},$ where $Y_k \sim N(\mu_{k}, \sigma^2)$  with unknown mean $\mu_k$ and known variance $\sigma^2$ for $k=1, \ldots,K.$
Suppose we are interested in estimating the $K$-dimensional mean vector $\mu= (\mu_1, \ldots, \mu_{K})^{T}$ and want to find a good estimator $\hat \mu = (\hat \mu_1, \cdots, \hat \mu_{K})^{T}$ under the mean squared error (MSE) criterion $MSE(\hat \mu) = \E ||\hat \mu - \mu||^2 =  \E \left(\sum_{k=1}^{K} (\hat \mu_{k} - \mu_{k})^2 \right).$


It is trivial to see that the method of moment estimator (MOM) or maximum likelihood estimator (MLE) of $\mu_{k}$ is $\hat{\mu}_k^{MLE} = Y_k$ for $k=1, \ldots, K,$ since each $\mu_{k}$ corresponds to only one normal variable $Y_{k}.$ A surprising result in a remarkable paper by James and Stein \cite{james1961estimation} is that there are uniformly better estimators than MOM or MLE in the sense of smaller MSE when simultaneously estimating $K \ge 3$ unknown parameters! Since then shrinkage estimation has become a basic tool in the analysis of high-dimensional data, especially when the object to estimate holds sparsity properties.

Many kinds of shrinkage estimators have been developed in the literature, see Cand\'{e}s \cite{candes2006modern} for the review and more references. Below we will review two kinds of shrinkage estimators that will be used in our proposed quickest change detection schemes. The first one is the linear shrinkage estimator of $\mu_k$'s defined by
\begin{eqnarray} \label{eq:mujs}
\hat{\mu}_k &=& a Y_{k} + (1-a) \zeta \cr
&=& a Y_{k} + b,
\end{eqnarray}
where $0 \le a \le 1$ is the shrinkage factor, $\zeta$ is a pre-specified real-valued constant (e.g., $\zeta= 0$), and $b = (1-a) \zeta.$ This corresponds to shrinking the observed vector $(Y_1,\cdots, Y_K)^{T}$ to the pre-specified vector $(\zeta, \cdots, \zeta)^{T}$ as the shrinkage factor $a$ goes to $0$ (note that in a more general setting,  $\zeta$ can be different for different $k$). 
Observe that the linear shrinkage estimator $\hat{\mu}_k$  in (\ref{eq:mujs}) has the common shrinkage factor $a$ for all $k,$ and intuitively this works well when all true $\mu_{k}$'s are nonzero, or better yet, have  similar values. The second kind of shrinkage estimator  is the hard-thresholding estimator defined by
\begin{eqnarray} \label{for: thresh01}
  \hat{\mu}_k  = \left\{
\begin{array}{ll}
   Y_k  & \hbox{ if $|Y_k| \geq \omega_k$} \\
   \mu_0 = 0 & \hbox{ if $ |Y_k| < \omega_k$}
\end{array}
\right. .
\end{eqnarray}
Intuitively, the  hard-thresholding estimator in (\ref{for: thresh01})  works when only a (small) subset of $\mu_k$'s are different from $\mu_0 = 0$. In such scenario, it makes more sense to shrinking non-significant MOM or MLE estimators of $\mu_k$'s directly to $0.$ Indeed, the optimality properties of hard-thresholding estimators in (\ref{for: thresh01}) were established in the context of offline point estimation, see, for example, Donoho and Johnstone \cite{donoho1994ideal}.


\par

\section{Our proposed monitoring schemes}

In the problem of online monitoring of $K$ independent normally distributed data streams with possible mean changes, if we completely knew each and every post-change parameter $\mu_{k},$ then many classical quickest change detection procedures for monitoring one-dimensional data stream can be easily adapted to develop global monitoring schemes, and one of them is the well-known Shiryaev-Roberts procedure  (Shiryaev \cite{shiryaev1963optimum} and Roberts \cite{roberts1966comparison}) that can be defined as follows in our context.  Let $\Lambda^{SR}_{n,m}$ be the likelihood ratio statistic of all observations up to time $n$ in the problem of testing $H_0:$ {\it no change} against $H_1:$ {\it a change occurs at time $m (\le n),$} i.e.,
\begin{eqnarray} \label{eqnSec03e01}
 \Lambda^{SR}_{n,m} = \prod_{ \ell = m}^{n}  \prod_{k = 1}^{K} \frac{f_{\mu_{k}} (X_{k, \ell})}{f_{\mu_0} (X_{k,\ell})},
\end{eqnarray}
where $f_{\mu}(\cdot)$ is the probability density function of $N(\mu, 1).$ At time $n,$ the Shiryaev-Roberts procedure computes the global monitoring statistics
\begin{eqnarray} \label{SR_stat}
R^{SR}_{n} = \sum_{m = 1}^{n} \Lambda^{SR}_{n,m},
\end{eqnarray}
which can be thought of as assigning a uniform prior on the potential change-point values $\nu = m \in \{1,2, \cdots, n\}.$  Then
 the Shiryaev-Roberts procedure raises a global alarm at time
\begin{eqnarray} \label{for: change rule}
N_B^{SR} = \inf \{n \ge 1:  R^{SR}_{n} \geq B\},
\end{eqnarray}
where the threshold $B > 0$ is chosen to satisfy the ARL to false alarm constraint in (\ref{eq: change detect cond}).

When the post-change parameters $\mu_{k}$'s are unknown, one natural possibility is
to replace them by their corresponding estimators from the observed data. In the quickest change detection literature, it is standard to use MLE or MOM to estimate the unknown post-change parameters, though there are generally two different approaches, depending on whether or not to use the same estimate for all $n-m+1$ post-change parameters $\mu_{k}$'s for $\ell = m, m+1, \ldots, n$  in the likelihood ratio $\Lambda^{SR}_{n,m}$ in (\ref{eqnSec03e01}) at time $n (\ge m).$
The first one is to replace all $n-m+1$ $\mu_{k}$'s by the same estimator based on all observations from the putative change-point time $m$ to the current time step $n$, and thus it often leads to the generalized likelihood ratio type statistic, see Xie and Siegmund \cite{xie2013sequential}.

The second approach, adopted  by Lorden and Pollak \cite{lorden2005nonanticipating},  is to use different estimates to the $n-m+1$ $\mu_{k}$'s. To be more concrete, for each $k=1, \ldots, K,$ Lorden and Pollak \cite{lorden2005nonanticipating}  considers $n-m+1$ MLE/MOM estimates of $\mu_{k}:$
\begin{eqnarray} \label{SRRS_est}
\hat{\mu}_{k,m,\ell} &=& {\bar  X}_{k,m,\ell} \\
& = &
\left\{
  \begin{array}{ll}
   \frac{X_{k,m}+ ...+ X_{k,\ell-1}}{\ell-m}, & \hbox{if $\ell= m+1, \ldots, n$} \\
  \mu_0 = 0, & \hbox{if $\ell = m$}
  \end{array}
\right.  \nonumber
\end{eqnarray}
and then proposes to plug these $\hat{\mu}_{k,m,\ell}$'s into (\ref{eqnSec03e01})-(\ref{for: change rule}) to yield the quickest change detection scheme. It is important to note that at time $\ell,$ the estimate $\hat{\mu}_{k,m,\ell}=  {\bar  X}_{k,m,\ell}$ in (\ref{SRRS_est}) only uses the observations, $X_{k,m}, \ldots, X_{k, \ell-1},$ to estimate $\mu_{k}$ at time $\ell,$ which allows one to reserve the observation $X_{k, \ell}$ only for detection of a change. By doing so, we keep  two important properties of $\Lambda^{SR}_{n,m}$ in (\ref{eqnSec03e01}): (i) the recursive form $\Lambda^{SR}_{n,m} = \Lambda^{SR}_{n-1,m}  \prod_{k = 1}^{K} [f_{\mu_{k}} (X_{k, n})/f_{\mu_0} (X_{k,n})],$ and (ii) the nice property of $\E_{\infty}(\Lambda^{SR}_{n,m}) = 1$ which leads to a useful fact that $R^{SR}_{n} - n$ is a martingale under the pre-change hypothesis. Lorden and Pollak \cite{lorden2005nonanticipating} termed their scheme as Shiryaev-Roberts-Robbins-Siegmund (SRRS) scheme, as similar idea  has been used earlier in Robbins and Siegmund \cite{robbins1972class} for sequential hypothesis testing problems. Below the scheme of  Lorden and Pollak \cite{lorden2005nonanticipating} will be called as the original SRRS scheme, and will be denoted by $N^{orig}_{B}.$ It was shown in Lorden and Pollak \cite{lorden2005nonanticipating}  that the original SRRS scheme $N^{orig}_{B}$ is first-order asymptotically optimal when monitoring $K=1$ data stream as the ARL to false alarm constraint $A$ in (\ref{eq: change detect cond}) goes to $\infty.$ After a careful analysis, it can also be shown that the first-order asymptotic optimality properties of  $N^{orig}_{B}$ can be extended for any fixed dimension $K,$ but unfortunately
the second-order term of the detection delay of the original SRRS scheme $N^{orig}_{B}$  is a linear function of $K.$  In other words, the original SRRS scheme $N^{orig}_{B}$ of Lorden and Pollak \cite{lorden2005nonanticipating} suffers the same problem of many classical schemes mentioned in Mei \cite{mei2010efficient} that the coefficient of the second-order term of detection delay is of order $K,$ and  thus its first-order asymptotic  optimality properties can be meaningless in the practical setting of monitoring large-scale data streams.

In this paper, we propose to develop a global monitoring scheme by combining the shrinkage estimators with the SRRS scheme of Lorden and Pollak \cite{lorden2005nonanticipating}. Our motivation is fueled by the fact that we need to estimate $K$ post-change means $\mu_{k}$'s simultaneously: if we let $Y_{k} = \hat{\mu}_{k,m,\ell}$ in (\ref{SRRS_est}) for all $k=1, \ldots, K,$ then existing research in the offline point estimation suggests that the shrinkage estimators in (\ref{eq:mujs}) or (\ref{for: thresh01}) should lead a better estimation of the true unknown post-change means $\mu_{k}$'s, which might lead to a better quickest change detection scheme.

Inspired by the linear shrinkage estimator in (\ref{eq:mujs}) and the hard thresholding estimator in (\ref{for: thresh01}), we propose a systematic approach that performs the linear shrinkage for values of MLE/MOM ${\bar X}_{k,m,\ell}$'s in (\ref{SRRS_est}) that are not thresholded. Specifically, 
we propose to  consider the shrinkage estimators of the form
\begin{eqnarray} \label{eq: thresh def}
  \hat{\mu}_{k,m,\ell} = \left\{
\begin{array}{ll}
   a {\bar X}_{k,m,\ell} + b  & \mbox{if $\ell=m+1, \cdots, n,$ and } \\
    &       \hbox{ $|{\bar X}_{k,m,\ell} | \geq \omega_k$}  \\
   c & \hbox{otherwise.}
\end{array}
\right. ,
\end{eqnarray}
where $a, b, c$ are three constants to be specified later. Note that $a=1, b =0$ and $c= 0$ correspond to the hard-thresholding estimators in (\ref{for: thresh01}), which will be shown later to be one of reasonable good choices under the post-change hypothesis $\Omega$ in  (\ref{par:par}).

Our proposed shrinkage-based SRRS schemes are defined by plugging the shrinkage/thresholding estimators $\hat{\mu}_{k,m,\ell}$ in (\ref{eq: thresh def}) into (\ref{eqnSec03e01})-(\ref{for: change rule}).  To be more concrete, define
\begin{eqnarray}  \label{new_LRaa}
\Lambda_{n,m} &=&  \prod_{ \ell = m}^{n} \prod_{k = 1}^{p} \frac{f_{\hat{\mu}_{k,m,\ell}} (X_{k,\ell})}{f_{\mu_0} (X_{k,\ell})} \\
&=& \Lambda_{n-1,m} \prod_{k = 1}^{p} \frac{f_{\hat{\mu}_{k,m,n}} (X_{k,n})}{f_{\mu_0} (X_{k,n})} \mbox{ for } n > m, \nonumber
\end{eqnarray}
where $\Lambda_{n,n}=1$ for all $n=1, 2, \ldots,$ and
\begin{eqnarray}
R_{n} &=& \sum_{m=1}^{n} \Lambda_{n,m},  \label{new_LR}
\end{eqnarray}
with $R_1 = 1.$ Then our proposed  shrinkage-based SRRS scheme  raises an alarm at the first time
\begin{eqnarray}\label{eqnSec03e07}
N_B = \inf \{n\ge 1: R_n \geq B\},
\end{eqnarray}
where $B > 0$ is a pre-specified threshold.

Note that   the original SRRS scheme in Lorden and Pollak \cite{lorden2005nonanticipating} can be thought of as a limiting case of our proposed  shrinkage-based scheme $N_B$ in (\ref{eqnSec03e07}) when $a=1, b= 0$ and $\omega_{k} \rightarrow 0$ in (\ref{eq: thresh def}). In addition, many arguments in the asymptotic analysis of the original SRRS scheme in Lorden and Pollak \cite{lorden2005nonanticipating} for $K=1$ dimension such as martingale properties and non-linear renewal theory for overshoot analysis can be applied to the proposed shrinkage-based scheme $N_B,$ subject to a careful analysis of the shrinkage estimators  in  (\ref{eq: thresh def}). Our major contribution is to introduce the shrinkage estimators to the quickest change detection problem and demonstrate its usefulness to lessen the dimension effects when the number $K$ of data streams is large.

%

\section{Asymptotic properties}

In this section, we investigate the asymptotic properties of the proposed shrinkage-based SRRS scheme $N_B$ in  (\ref{new_LR})  and  (\ref{eqnSec03e07}) when the estimators   $\hat{\mu}_{k,m,\ell}$'s of the post-change means $\mu_{k}$'s are the shrinkage estimators in (\ref{eq: thresh def}). The following discussion is divided into three subsections. The first two subsections address two properties of the proposed shrinkage-based SRRS scheme under the general setting: the ARL to false alarm and detection delay, respectively. The third subsection focuses on the suitable choice of tuning parameters in our proposed shrinkage-based SRRS scheme.



\subsection{The ARL to false alarm}

To derive the ARL to false alarm of the proposed shrinkage-based SRRS scheme $N_B$ in (\ref{new_LR}) and (\ref{eqnSec03e07}), it is crucial to  observe that its global monitoring statistic $R_n$ is the Shiryaev-Roberts-type statistics and thus $R_n - n$ is a martingale under the pre-change hypothesis. By the well-known Doob's optional stopping time theorem (see Theorem 10.10 of Williams \cite{williams1991probability}), for the stopping time $N= N_B$ defined in (\ref{eqnSec03e07}), we have $\E_{\infty}(N) = \E_{\infty}(R_N) \ge B,$ as $R_{N_B} \ge B$ by the definition of $N_B.$ Also see the proof of Theorem 4 of Lorden and Pollak \cite{lorden2005nonanticipating} for more detailed arguments.  The following theorem summarizes this result.
\begin{theorem} \label{thm1}
Consider the proposed  shrinkage-based SRRS scheme $N_B$ in (\ref{new_LR}) and (\ref{eqnSec03e07}) with $\hat{\mu}_{k,m,\ell}$ being the shrinkage estimators in (\ref{eq: thresh def}). For any $B > 0,$
\[
\E_{\infty}(N_B) \ge B.
\]
\end{theorem}

\medskip
While Theorem \ref{thm1} is applicable regardless of the value of the dimension $K$ (the number of data streams), it is important to point out that the Monte Carlo simulation of $\E_{\infty}(N_B)$ is a different story due to the curse of dimensionality. If the dimension $K$ is small, say $K=1$ or $5,$ then a Monte Carlo simulation with runs of thousands will provide a reasonable estimate of $\E_{\infty}(N_B)$ for a moderately large threshold $B,$ say $B= 10^4.$ However, the number of necessary runs is exponentially increasing as the dimension $K$ increases, as the scheme $N_B$ is highly skewed for large $K,$ and the sample mean or median based on $10^5$ or $10^6$ of realizations of $N_B$ can be a very poor estimate of $\E_{\infty}(N_B).$

The reason is that the likelihood ratio $\Lambda_{n,m} (m < n)$  in (\ref{new_LRaa}) and the global monitoring statistic $R_{n}$ in (\ref{new_LR}) are typically highly skewed to $0$ and $1$ for large dimensional $K,$ respectively. To see this, consider the likelihood ratio $\Lambda_{n,m}$ when $\hat{\mu}_{k,m,\ell}$ is the MLE/MOM estimates in (\ref{SRRS_est}). On the one hand, for a fixed $n$ and any given $1 \le m \le n-1,$ we have $\E_{\infty}( \Lambda_{n,m} ) = 1$ and $\E_{\infty}(R_n) = n.$ On the other hand, for normal distributions,
\begin{eqnarray*}
& & \E_{\infty}\log(\Lambda_{n,m} ) \\
&=& \sum_{\ell=m}^{n} \sum_{k=1}^{K} \E_{\infty}   \Big( \E_{\infty}(\hat \mu_{k,m,\ell} X_{k, \ell} - \frac12 (\hat \mu_{k,m,\ell})^2 \Big| \hat \mu_{k,m,\ell}  \Big) \\
&=& -\frac12 \sum_{\ell=m}^{n} \sum_{k=1}^{K} \E_{\infty} (\hat \mu_{k,m,\ell})^2  \quad \mbox{ (as $\E_{\infty}(X_{k, \ell}) = 0$) }  \\
&=& - \frac12 \sum_{k=1}^{K} \sum_{\ell=m+1}^{n} \frac{1}{\ell - m}  \quad \mbox{ (as $\hat \mu_{k,m,m} = 0$) } \\
&=& - \frac{1}{2}  K  (1+ \frac{1}{2} + \cdots + \frac{1}{n-m}).
\end{eqnarray*}
Here the third equation uses the fact that when $\hat{\mu}_{k,m,\ell}$ is the MLE/MOM estimates in (\ref{SRRS_est}), it has a $N(0, 1/(\ell-m))$ distribution. Now when $K=100,$ we have $\E_{\infty}\log(\Lambda_{n,n-1} ) = - 50,$ implying that $\Lambda_{n,n-1}$  is concentrated around $e^{-50} = 1.9 \times 10^{-22},$ even though  $\E_{\infty}(\Lambda_{n,n-1}) = 1.$ For all other $m < n,$ the likelihood ratios $\Lambda_{n,m}$'s will be concentrated around an even smaller value. Hence, for a fixed time $n,$ when we simulate the global monitoring statistic $R_{n}$ of the original SRRS scheme,  we will mostly likely observe $R_{n} \approx 1$  (recall that $\Lambda_{n,n}$ is defined as a constant $1$), although $\E_{\infty}(R_n) = n.$

The above argument can also be extended to the proposed SRRS scheme with the shrinkage estimators in (\ref{eq: thresh def}), and our numerical experiences seem to suggest that the Monte Carlo estimate of $\E_{\infty}(N_B)$ works poorly and is highly biased unless the linear shrinkage factor $a$ is of order $O(1/\sqrt{K}).$ As mentioned in Rubinstein and Glynn \cite{rubin}, the curse of dimensionality is one of the central topics in Monte Carlo simulation due to the degeneracy properties of likelihood ratios, and the importance sampling technique does not help in the high dimensional problem 
unless we can reduce it to an equivalent low-dimension problem. It remains an open problem how to overcome  the curse of dimensionality to simulate $\E_{\infty}(N_B)$ effectively for our proposed SRRS scheme $N_{B}$  in the general context of monitoring a large number $K$ of data streams.

A challenging practical question is how to find the threshold $B$ of the proposed shrinkage-based SRRS scheme $N_B$ in (\ref{new_LR}) and (\ref{eqnSec03e07}), so that it satisfies the pre-specified ARL to false alarm constraint $A$ in (\ref{eq: change detect cond}). The good news is that Theorem 1 provides a theoretical bound: a choice of $B=A$ will guarantee that the proposed shrinkage-based SRRS scheme $N_B$ satisfies the ARL to false alarm constraint in (\ref{eq: change detect cond}). For that reason,  in our numerical simulations below, we will set $B = A$ and report the impact of shrinkage estimation on the detection delays of the proposed shrinkage-based SRRS scheme $N_B.$


\medskip
\subsection{Detection delay}

In this subsection, we derive the asymptotic expression of the detection delay of the proposed shrinkage-based SRRS scheme $N_B$ in (\ref{new_LR}) and (\ref{eqnSec03e07}) under the setting when the dimension $K$ is fixed and the threshold $B$ goes to $\infty.$ In this subsection and only in this subsection, we assume that a change occurs to the $k$-th data stream at time $\nu$ and the true post-change mean of the $k$-th data stream is $\mu_{k}$ for all $k=1, \ldots, K.$ That is, the true post-change mean vector $\mu=(\mu_1, \ldots, \mu_{K})^{T}.$ Recall that if $\mu_{k} = \mu_0$ then no changes occur to the $k$-th data streams.

To present the results on the detection delays, we need to first introduce some new notation.
For each $k=1, \ldots, K,$  denote by $\mu_k^{*}$ the limit of  the shrinkage estimators $\hat{\mu}_{k,m,\ell}$ in (\ref{eq: thresh def}) as $\ell \rightarrow \infty,$ i.e.,
$\mu_k^{*} = \lim_{\ell \rightarrow \infty} \hat{\mu}_{k,m,\ell}$ under the post-change hypothesis. Note that the limit $\mu_k^{*}$ does not depend on the initial time $m$ of  the estimators, and for the shrinkage estimator $\hat{\mu}_{k,m,\ell}$ in (\ref{eq: thresh def}), it is easy to see that
\begin{eqnarray} \label{eqn_0015}
\mu_{k}^{*} =  \left\{
\begin{array}{ll}
   a \mu_{k} + b & \hbox{if $|\mu_{k}| > \omega_k$} \\
    c & \hbox{if $|\mu_{k}| < \omega_{k}$}
\end{array}
\right. ,
\end{eqnarray}
for each $k=1, \ldots, K.$ Here we purposely do not consider the cases of $\mu_{k} = \pm \omega_{k},$ as the corresponding analysis is complicated since the corresponding limit $\mu_{k}^{*}$  can be either $a \mu_{k} + b$ or $c,$ either with probability $1/2.$
This is also the reason why the post-change hypothesis set $\Omega$ in (\ref{par:par}) makes a technical assumption that the post-change mean $\mu_{k} \ne \pm \omega_{k}$ so as to simplify our theoretical analysis.

Denote the vector of the limits $\mu_{k}^{*}$'s in (\ref{eqn_0015}) by  $\mu^* = (\mu_{1}^{*}, \ldots, \mu_{K}^{*})^{T}.$ It is important to note that in Lorden and Pollak \cite{lorden2005nonanticipating}, or more generally in the quickest change detection literature, the limit vector $\mu^{*}$ is always the same as the true post-change mean vector $\mu$. Hence,  the asymptotic analysis on the detection delay of the scheme $N_B$ in (\ref{new_LR}) and (\ref{eqnSec03e07}) is closely related to the classical Shiryaev-Roberts procedure $N_{B}^{SR}$ that detects a change from $\mu_0$  to the known post-change $\mu.$
However, for our proposed shrinkage estimators in  (\ref{eq: thresh def}), it is no longer true that  $\mu^{*} = \mu$. Hence, we need to compare $N_B$ with the Shiryaev-Roberts procedure that mis-specifies the post-change means, i.e., the one that is designed to detect a change from $\mu_0$  to $\mu^{*}$ but the true post-change mean vector is actually $\mu.$

For that reason,  we define a new information number
\begin{eqnarray} \label{eqn06sec3-info1}
I(\mu^{*}, \mu_0; \mu) &=&  \E_{\mu} \sum_{k=1}^{p} \left ( \log \frac{f_{\mu_{k}^{*}}(X_{k,\ell})}{f_{\mu_0}(X_{k,\ell})} \right ).
\end{eqnarray}
When $f_{\mu_0}$ and $f_{\mu_{k}}$ are normal distributions with common variance $1,$ it becomes
\begin{eqnarray} \label{eqn06sec3-info}
I(\mu^{*}, \mu_0; \mu) &=& - \frac12 \sum_{k=1}^{p} \mu_{k}^{*}  (\mu_{k}^{*} - 2 \mu_{k}).
\end{eqnarray}
Plugging the limits $\mu_{k}^{*}$'s in (\ref{eqn_0015}) directly into (\ref{eqn06sec3-info}) yields that
\begin{eqnarray} \label{eqn02sec5-info_gen}
I(\mu^{*}, \mu_0; \mu) &=& -\frac12 \sum_{k: |\mu_{k}| >  \omega_k} \big(a  \mu_k + b \big) \big( (a - 2) \mu_{k}  + b \big) \cr
&& - \frac12 \sum_{k: |\mu_{k}| <  \omega_k} c(c - 2 \mu_{k}).
\end{eqnarray}
In the special case when $(a,b,c) =(1,0,0),$ the new information number $I(\mu^{*}, \mu_0; \mu)$ in (\ref{eqn02sec5-info_gen}) has a simpler form:
\begin{eqnarray} \label{eqn02sec5-info}
I(\mu^{*}, \mu_0; \mu) &=&   \frac12 \sum_{k: |\mu_{k}| >  \omega_k} (\mu_k)^2.
\end{eqnarray}
If we  further let $\omega_{k}$ go to $0,$ this becomes a more familiar form
\begin{eqnarray} \label{eq:Itot}
I_{tot} = \frac12 \sum_{k=1}^{K} (\mu_{k})^2,
\end{eqnarray}
where we change the notation to $I_{tot}$ to emphasize its special meaning as the  Kullback-Leibler divergence between the pre-change and post-change hypotheses. In Information Theory and Statistics, the  Kullback-Leibler divergence $I_{tot}$ in (\ref{eq:Itot})  has been regarded as a measure to characterize the distance between pre-change and post-change distributions, or equivalently, as a measure how difficulty it is to detect the change. However, in the context of monitoring a large $K$ number of data streams,  $I_{tot}$ in (\ref{eq:Itot}) might  no longer be as informative as one thought, since it ignores the spatial uncertainty associated with which subset of data streams are affected. A more meaningful measure that takes into account such a spatial uncertainty will be $I(\mu^{*}, \mu_0; \mu)$ in (\ref{eqn02sec5-info}), or more generally those in (\ref{eqn02sec5-info_gen}) or (\ref{eqn06sec3-info1}), which can be thought of as the shrinkage version of  the  Kullback-Leibler divergence $I_{tot}$ in (\ref{eq:Itot}).

With the notation in (\ref{eqn02sec5-info_gen})-(\ref{eq:Itot}), we are ready to present our main result on the detection delays of the proposed shrinkage-based SRRS scheme $N_B$ in (\ref{new_LR}) and (\ref{eqnSec03e07}).
\begin{theorem} \label{thm: detection change}
Consider the proposed shrinkage-based SRRS scheme $N_B$ in (\ref{new_LR}) and (\ref{eqnSec03e07}) with the estimators $\hat{\mu}_{k,m,\ell}$'s defined in (\ref{eq: thresh def}). Assume that the information number $I(\mu^{*}, \mu_0; \mu)$ in (\ref{eqn02sec5-info_gen}) is positive. Then, as $B \rightarrow \infty,$ the detection delay of $N_{B}$ satisfies
\begin{eqnarray}  \label{eqn05sec4}
   D_{\mu} (N_B)  &\le& \frac{\log B}{I(\mu^{*}, \mu_0; \mu)} + \frac{r a^2}{2 I(\mu^{*}, \mu_0; \mu)} \times \cr
  && \times \log\Big( \frac{\log B}{I(\mu^{*}, \mu_0; \mu)} \Big) + o( r \log \log B),\quad
\end{eqnarray}
where $r$ is defined in (\ref{par:par3}) and represents the true number of affected data streams with ``big" local changes.
\end{theorem}

\medskip
The detailed proof of this theorem will be presented in the Appendix.   Note that the original SRRS scheme $N_B^{orig}$  proposed in Lorden and Pollak \cite{lorden2005nonanticipating}  can be thought of as the special case of our proposed scheme $N_B$ in (\ref{new_LR}) and (\ref{eqnSec03e07}). 
By Theorems 1 and 2, or by extending the proof of Theorem 4 in Lorden and Pollak \cite{lorden2005nonanticipating} from one-dimensional to $K$-dimensional, we can establish the first-order asymptotic optimality of the original SRRS scheme $N_{B}^{orig}$  in Lorden and Pollak \cite{lorden2005nonanticipating} as follows:

\begin{corollary} \label{cor1}
 As $B \rightarrow \infty,$ the original SRRS scheme $N_B^{orig}$ satisfies
\begin{eqnarray}  \label{eqn06sec4}
    D_{\mu} (N_B^{orig})  &\le&  \frac{\log B}{I_{tot}} +  \frac{K}{2 I_{tot}} \log\Big( \frac{\log B}{I_{tot}} \Big)  \cr
   && + o(K \log \log B),
\end{eqnarray}
where  $I_{tot}$ is defined in (\ref{eq:Itot}). Moreover, if we let $B=A,$ then the original SRRS scheme $N_B^{orig}$  is first-order asymptotically optimal in the sense of asymptotically minimizing the detection delay  $D_{\mu}(N_B^{orig})$ for each and every post-change mean vector $\mu$ subject to the false alarm constraint in (\ref{eq: change detect cond})  when $K$ is fixed and the constraint $A$  in (\ref{eq: change detect cond}) goes to $\infty.$
\end{corollary}

\proof For  the original SRRS scheme $N_B^{orig},$ the limit $\mu_{k}^{*}$ in (\ref{eqn_0015}) becomes the true post-change mean $\mu_{k}$ itself, and thus it is clear that  the original SRRS scheme $N_B^{orig}$  can be thought of as the special case of our proposed scheme $N_B$ in (\ref{new_LR}) and (\ref{eqnSec03e07}) with $a=1, r=K$ and $I(\mu^{*}, \mu_0; \mu) = I_{tot}.$ Relation (\ref{eqn06sec4}) then follows directly from Theorem 2. 
By Theorem 1, the choice of $B=A$ makes sure that the original SRRS scheme $N_B^{orig}$ satisfies the false alarm constraint in (\ref{eq: change detect cond}). Then the asymptotic optimality properties of  $N_B^{orig}$  follows at once from a well-known lower bound  on the detection delay of any scheme $N$ satisfying the false alarm constraint in (\ref{eq: change detect cond}):  $D_{\mu} (N) \ge (1+o(1))(\log A) / I_{tot},$ see  Lorden \cite{lorden1971procedures}.
\endproof

\medskip
While Corollary 1 establishes the first-order asymptotic optimality property of the original SRRS scheme $N_{B}^{orig},$ 
it can be meaningless in practical setting when the dimension $K$ is large and $B$ is only moderately large. This is because
the second-order term ($\log\log(B)$) in the right-hand side of (\ref{eqn06sec4}) has coefficient $K$ and can be significant as compared to the first-order term $\log(B).$ A comparison of  (\ref{eqn05sec4}) and (\ref{eqn06sec4}) shows that shrinkage estimators impact the detection delays in two different places: one is the information number $I(\mu^{*}, \mu_0; \mu)$  in (\ref{eqn02sec5-info_gen}) on the first-order term, and the other is the factor $r a^2$ in the second-order term. These will allow us to illustrate in the next subsection how a suitable choice of shrinkage estimators in (\ref{eq: thresh def}) can reduce the overall detection delay.

\subsection{How to choose suitable shrinkage estimators?}

In our proposed shrinkage-based SRRS scheme $N_B$ in (\ref{new_LR}) and (\ref{eqnSec03e07}) with the estimators $\hat{\mu}_{k,m,\ell}$'s defined in (\ref{eq: thresh def}), there are two sets of tuning parameters: one is the lower bounds $\omega_{k}$'s and the other is the constant $(a,b,c)$. The choices of the lower bounds $\omega_{k}$'s are straightforward, as they are pre-specified in the post-change hypothesis set  $\Omega$ in (\ref{par:par}). Below we will focus on the suitable choice of tuning parameter $(a,b,c).$

By Theorem 2, if we want to minimize the first-order term of the detection delay of the proposed shrinkage-based scheme $N_B,$
then it suffices to maximize  the information number $I(\mu^{*}, \mu_0; \mu)$ in (\ref{eqn02sec5-info_gen}). Hence, it is natural to define the ``first-order" optimal choice of $(a, b , c)$ as the one that maximizes $I(\mu^{*}, \mu_0; \mu)$ in (\ref{eqn02sec5-info_gen}). The following theorem provides the corresponding ``first-order" optimal choice of $(a,b,c)$ among all possible shrinkage estimators $\hat{\mu}_{k,m,\ell}$'s in (\ref{eq: thresh def}):

\medskip
\begin{theorem} \label{cor2b}
Under the post-change hypothesis set $\Omega$ in (\ref{par:par}), the choice of $a=1, b=0, c=0$ is ``first-order" optimal for the proposed SRRS scheme among all possible shrinkage estimators $\hat{\mu}_{k,m,\ell}$'s in (\ref{eq: thresh def}).
\end{theorem}

\proof It suffices to show that $a=1, b=0, c=0$ maximizes $I(\mu^{*}, \mu_0; \mu)$ in (\ref{eqn02sec5-info_gen}). Note that the right-hand side of (\ref{eqn02sec5-info_gen}) is a quadratic function of $a,b,c$ and thus the optimal values can be found by taking derivatives of the right-hand side of  (\ref{eqn02sec5-info_gen}) with respect to $a,b,c.$ Following the definition of $r$ in (\ref{par:par3}), define
\begin{eqnarray*}
D_1 &=&  \sum_{k=1}^{K} \mu_{k} 1\{|\mu_{k}| >  \omega_k \} \\
D_2 &=& \sum_{k=1}^{K} \mu_{k}^2 1\{|\mu_{k}| >  \omega_k \}  \\
D_3 &=& \sum_{k=1}^{K} \mu_{k} 1\{|\mu_{k}| <  \omega_k \}.
\end{eqnarray*}
Then the derivatives of the right-hand side of  (\ref{eqn02sec5-info_gen}) with respect to $a,b,c$ can be rewritten as
\begin{eqnarray*}
 D_1 (a-1)  + r b &=& 0; \\
 D_2 (a-1)  + D_1 b  &=& 0; \\
 (K-r) c  - D_3 &=& 0.
\end{eqnarray*}
Clearly, under the post-change hypothesis $\Omega$ in (\ref{par:par}), the post-change mean $\mu_{k} = 0$ whenever $|\mu_{k}| \le \omega_{k},$ implying that $D_3 = 0.$ Hence,
$(a^* ,  b^*, c^*) = (1, 0, 0)$ is the unique optimal choice of $(a,b,c)$ when $(D_1)^2 \ne r D_2$ and $r \ne K,$ and is one of infinitely many optimal solutions otherwise. Thus the theorem holds.
\endproof

\medskip
When $a=1, b=0, c=0,$ the shrinkage estimators $\hat{\mu}_{k,m,\ell}$'s in (\ref{eq: thresh def}) become the hard-thresholding estimators
\begin{eqnarray} \label{eq:hard_def}
  \hat{\mu}_{k,m,\ell} = \left\{
\begin{array}{ll}
   {\bar X}_{k,m, \ell} 
    & \hbox{if $\ell = m+1, \ldots, n,$ and} \\
    & \hbox{$|{\bar X}_{k,m, \ell}| \geq \omega_k$}  \\
   \mu_0=0 & \hbox{otherwise}
\end{array}
\right. ,
\end{eqnarray}
where ${\bar X}_{k,m, \ell}$'s are the MLE/MOM estimates of $\mu_{k}$ in (\ref{SRRS_est}).
Denote by   $N_B^{hard}$ the corresponding SRRS scheme $N_B$ in (\ref{new_LR}) and (\ref{eqnSec03e07}) when the estimators $\hat{\mu}_{k,m,\ell}$'s being the hard-thresholding estimators (\ref{eq:hard_def}).  The following corollary summarizes its first-order asymptotic optimality properties:

\begin{corollary} \label{cor2}
For any fixed dimension $K,$ the hard-thresholding scheme $N_{B}^{hard}$ with $B=A$ asymptotically minimizes the detection delay $D_{\mu} (N_{B}^{hard})$ (up to first-order) for each and every post-change mean vector $\mu \in \Omega$ in (\ref{par:par}) subject to the false alarm constraint $A$ in (\ref{eq: change detect cond}) as the constraint $A$ goes to $\infty.$
\end{corollary}

\proof By Theorem 2 and Corollary 1, it suffices to show that for the hard-thresholding estimators in (\ref{eq:hard_def}), $I(\mu^{*}, \mu_0; \mu)$ in (\ref{eqn02sec5-info}) is the same as $I_{tot}$ in (\ref{eq:Itot}) when $\mu \in \Omega$ in (\ref{par:par}). From the definition of $\Omega$ in (\ref{par:par}), we have $\mu_{k} = 0$ if $|\mu_{k}| < \omega_{k}.$ Thus $\sum_{k: |\mu_{k}| >  \omega_k} (\mu_k)^2 = \sum_{k=1}^{K} (\mu_{k})^2$  and it is clear from (\ref{eqn02sec5-info}) and (\ref{eq:Itot}) that $I(\mu^{*}, \mu_0; \mu) = I_{tot}$ for any $\mu \in \Omega.$ Hence the corollary holds.
\endproof

\medskip
It is useful to compare the original SRRS scheme $N_{B}^{orig}$ in Corollary 1 with the hard-thresholding scheme  $N_{B}^{hard}$ in Corollary 2.
On the one hand, the first-order asymptotic optimality property of  $N_{B}^{orig}$ is applicable to all possible post-change mean vectors $\mu$ no matter whether $\mu \in \Omega$ or not, whereas  $N_{B}^{hard}$ is first-order asymptotically optimal only for those $\mu \in \Omega$ in (\ref{par:par}). On the other hand, for these two schemes, the coefficients in the second-order terms of the detection delays are different: $K$ for $N_{B}^{orig},$ and  $r$ in (\ref{par:par3}) for $N_{B}^{hard}.$ This is exactly the  reason why the hard-thresholding estimators can reduce the detection delay in the sparse post-change case of $\Omega$ in (\ref{par:par})  when the number of affected data streams is much smaller than the total number of data streams, e.g., when $r=20$ out of $K=100$ data streams are affected. 

Corollary 2 also provides a partial answer to an open problem raised on page \#426 of Mei \cite{mei2010efficient} whether  we can develop new methods to reduce the coefficient in the second-order term of the detection delay from  $K$  to a smaller number while keeping the first-order asymptotic optimality properties. Our results show that such coefficient can be reduced to   the number  $r$ of affected data streams in the sparse post-change case.
We conjecture that $r$ in (\ref{par:par3}) is the smallest possible coefficient for the second-order term in the Gaussian model, but we do not have a rigorous proof.

Besides  the sparse post-change case, another interesting case of $\Omega$ in (\ref{par:par}) is when all data streams are affected simultaneously. In this case, we have $r = K,$ and thus the hard-thresholding scheme  $N_{B}^{hard}$  does not necessarily work efficiently, and to the best of our knowledge, no methodologies have been developed to improve the original SRRS scheme $N_{B}^{orig}$  or other classical quickest change detection schemes when the unknown local post-change means might be different for different local data streams.  Below we will demonstrate how to use Theorem 2 to derive a good choice of the linear shrinkage factor $a$ that can balance the tradeoff between the first-order and second-order of the detection delay.

To highlight our main ideas, let us focus on $a$ by setting $b=c= 0.$ Then the estimators $\hat{\mu}_{k,m,\ell}$'s  in (\ref{eq: thresh def}) becomes
\begin{eqnarray} \label{eq:linear_def}
  \hat{\mu}_{k,m,\ell} = \left\{
\begin{array}{ll}
   a {\bar X}_{k,m, \ell} 
    & \hbox{if $\ell = m+1, \ldots, n,$ and} \\
    & \hbox{$|{\bar X}_{k,m, \ell}| \geq \omega_k$}  \\
   \mu_0=0 & \hbox{otherwise}
\end{array}
\right. 
\end{eqnarray}
for some $0 \le a \le 1,$  where ${\bar X}_{k,m, \ell}$'s are the MLE/MOM estimates of $\mu_{k}$ in (\ref{SRRS_est}). Then $I(\mu^{*}, \mu_0; \mu)$ in (\ref{eqn02sec5-info_gen}) becomes
$I(\mu^{*}, \mu_0; \mu) = a(2-a) I_{tot}.$ By Theorem 2, when $r = K,$ minimizing the detection delay of $N_{B}$ is asymptotically equivalent to minimizing
\begin{eqnarray} \label{eq:maxa}
\frac{\log B}{a (2-a) I_{tot}} + \frac{K a}{2 (2-a) I_{tot}} \log\Big( \frac{\log B}{I_{tot}} \Big)
\end{eqnarray}
if we only keep the key terms containing the factor $a$ and ignore the $1/(a(2-a))$ factor inside the logarithm of the second term.  Clearly, $a=1$ maximizes $a(2-a),$ and this is equivalent to the first-order asymptotic optimality properties of $N_{B}^{orig}$ or $N_{B}^{hard}.$
However, a better choice of $a$ is to find $0 < a \le 1$ that minimize the summation in (\ref{eq:maxa}), not just the first term in (\ref{eq:maxa}). Note that a choice of $0 < a \le 1$ will make sure that the factor $a / (2-a)$ in the second term of (\ref{eq:maxa}) is less than  $1.$  The corresponding optimal value of $a$ will depend on $I_{tot},$ $\log B,$ and $K.$ For instance, when $B= 5000, K=100$ and $I_{tot} = 2.5,$ the summation in (\ref{eq:maxa}) becomes
\[
\frac{3.407}{a(2-a)} + \frac{24.516 a}{2-a}.
\]
This summation has the value $46.2$ when $a=1,$ and is minimized at $a=0.25$ with the smallest value $13.9.$ This suggests that a suitable choice of linear shrinkage estimators in (\ref{eq:linear_def})  can greatly reduce the overall detection delay as compared to the original SRRS scheme, although the price we pay is to sacrifice the first-order asymptotic optimality properties.


\par

\section{Numerical Simulations}

In this section, we report numerical simulations to illustrate the usefulness
of shrinkage or thresholding in the context of quickest change detection in Section V.A, and demonstrate the challenge of Monte Carlo simulations of the ARL to false alarms when  monitoring large-scale data streams in Section V.B.

\begin{figure}
  \centering
  \includegraphics[width=2.5in]{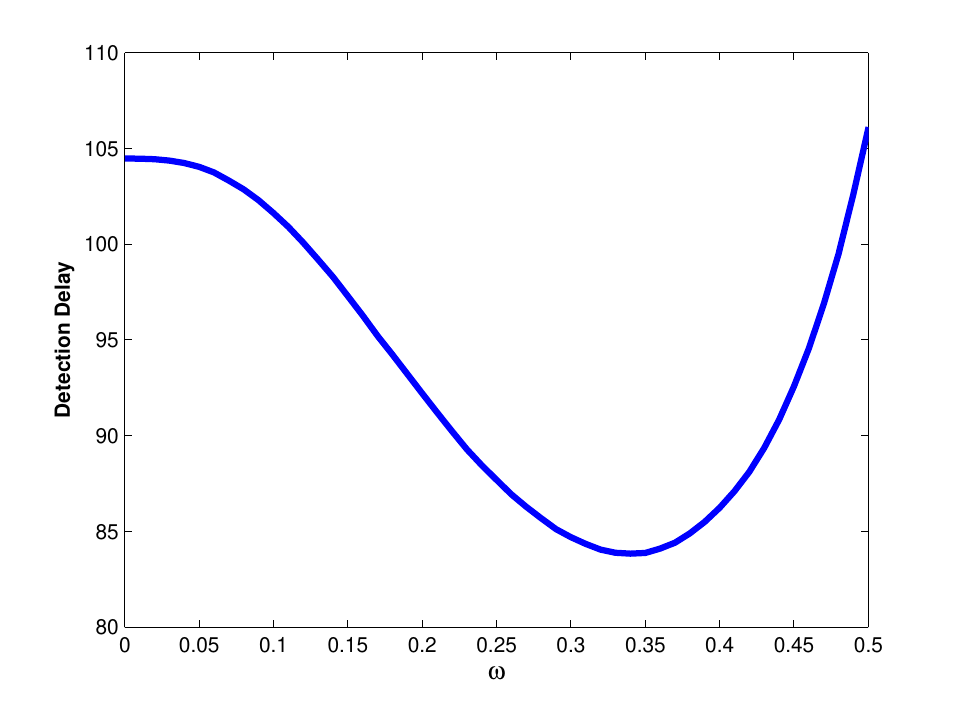}
\caption{The sparse post-change case with the hard-thresholding scheme $N_{B}^{hard}(\omega)$ that sets the common lower bound $\omega_{k} \equiv \omega$ for all $k=1, \ldots, K.$
The $x$-axis is the common lower bound $\omega$, and the $y$-axis is the simulated detection delay of $N_{B}^{hard}(\omega)$ when $B= 5000$. Here $\omega =0$ corresponds to the baseline scheme $N^{orig}_{B}$ without hard-thresholding. }

\label{figure:Itot}
\end{figure}

\begin{figure}
  \centering
 \includegraphics[width=2.5in]{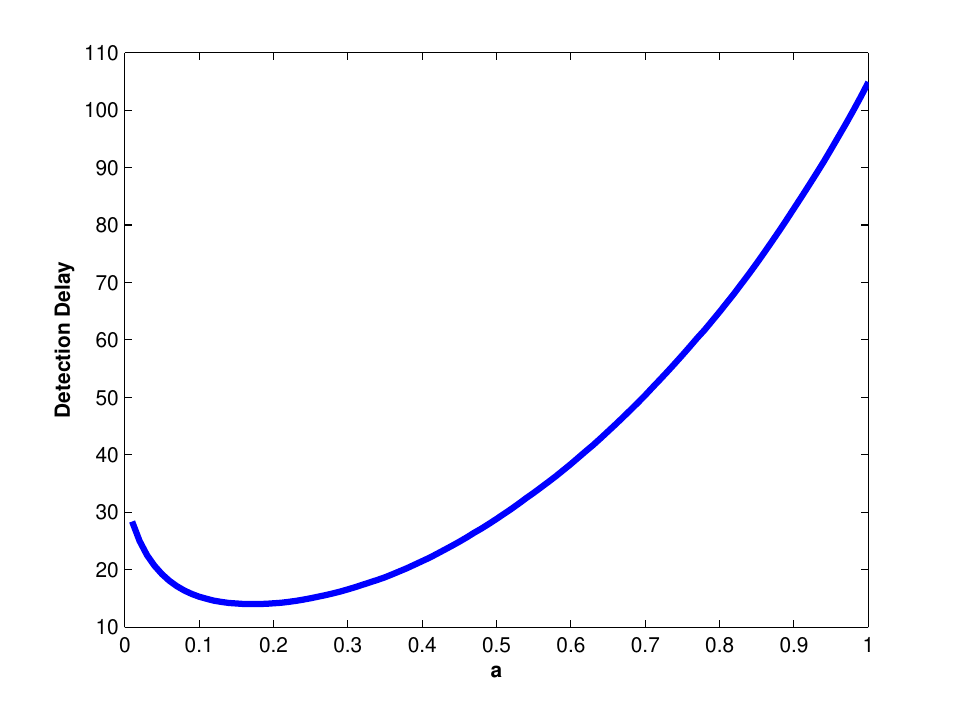} \\
  \includegraphics[width=2.5in]{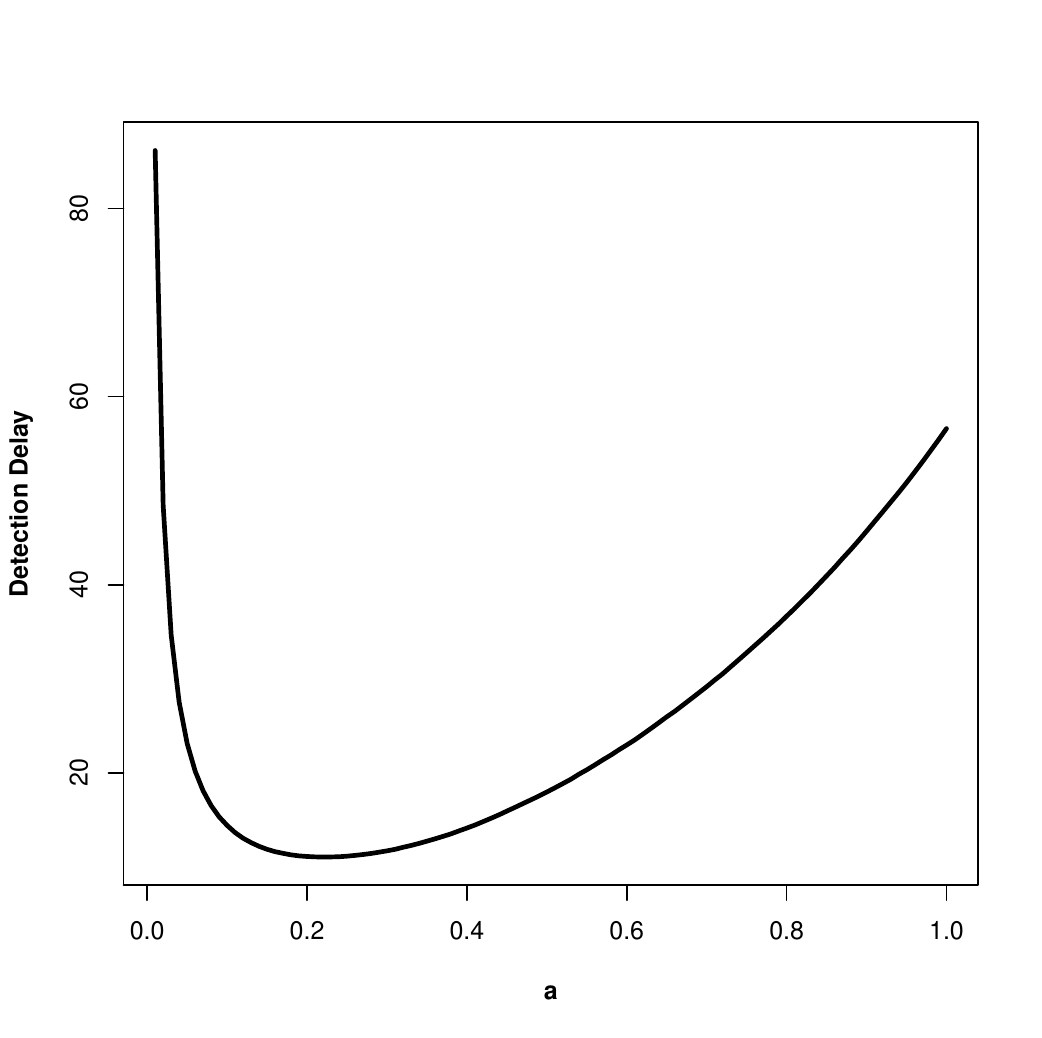}
\caption{The case when all data streams are affected, and we consider the SRRS scheme $N_{B}(a)$ with varied linear shrinkage factor $a$ while fixing $b=c=0$ with two different choices of fixed lower bounds $\omega_{k}$'s: the upper plot is when $\omega_{k} = 0$ for all $k,$ and the bottom plot is when $\omega_{k} = 0.01$ for all $k.$ The $x$-axis is the value of linear shrinkage factor $a$, and the $y$-axis is the simulated detection delay of $N_{B}(a)$ when $B= 5000$. Here $a=1$ corresponds to the scheme without linear shrinkage. }

\label{figure:Itotbb}
\end{figure}

\subsection{Shrinkage Effects}

Assume we are monitoring $K =100$ independent normal data streams whose initial distributions are $N(0,1)$ with possible changes in the means of some data streams. The ARL to false alarm constraint in (\ref{eq: change detect cond}) is assumed to be $\E_{\infty}(N) \ge A = 5000.$ As mentioned in subsection 4.1, we set all thresholds $B = 5000$ for all schemes $N_B$'s to avoid poor Monte Carlo estimates of $\E_{\infty}(N_B).$

We have conducted extensive simulations for different schemes under different kinds of post-change hypothesis $\Omega$ in   (\ref{par:par}), but will only report the results of two specific post-change hypotheses so as to highlight our findings.
The first one is the sparse post-change hypothesis case when $r=20$ out of $K=100$ data streams are affected, and the other is when all $K=100$ data streams are affected.  In both cases, we fix the overall information, $I_{tot}=\frac{1}{2}\sum_{k=1}^{K} \mu_k^2,$ to be $2.5.$ To be more specific, we consider two cases: (1)  when $r=20$ out of $K=100$ data streams are changed with the post-change mean $\mu_{k} = 0.5$ whereas there are no changes to the other remaining $K-r = 80$ data streams; and (2) when all $\mu_k = \sqrt{5/100} = 0.2236$ for all $k=1, \ldots, K.$ In both two cases, $I_{tot}=\frac{1}{2}\sum_{k=1}^{K} \mu_k^2=2.5.$
However, when we design the monitoring scheme, we will only know that the post-change mean vector $\mu$ is in (\ref{par:par}), and will  not use any other information of the true post-change parameters. As shown in the second remark on Page 1435 in Lorden and Pollak \cite{lorden2005nonanticipating}, the worst-case detection delays of the SRRS scheme $N_B$ occurs at time $\nu=1$, and thus we will report the detection delay performance of $N_{B=5000}$  under the post-change hypothesis when the change occurs at time $\nu = 1.$ All simulation results are based on $2500$ replications.

For the purpose of comparison, the baseline scheme is the original SRRS scheme $N_B^{orig}$ proposed by Lorden and Pollak \cite{lorden2005nonanticipating}. In the sparse post-change case when $r=20$ out of $K=100$ data streams are affected, we consider several different kinds of hard-thresholding schemes $N_{B}^{hard}$'s in Corollary 2. For the convenience of comparison, we set $\omega_k \equiv \omega$ for all $k,$ and then vary $\omega$ from $0$ (baseline) to $0.5$ with step size $0.01.$ For each hard-thresholding scheme with a given threshold $\omega,$ we then plot the detection delay of $N_{B}^{hard} = N_{B}^{hard}(\omega)$ as a function of $\omega$ in Figure \ref{figure:Itot}. It is evident from Figure \ref{figure:Itot} that the detection delay of the scheme $N_{B}^{hard}(\omega)$ is reduced from $104.9$ at $\omega = 0$ (baseline $N_B^{orig}$) to $83.8$ at $\omega = 0.35.$ This illustrates the usefulness of hard-thresholding estimators in the sparse post-change case.

In the case when all $K=100$ data streams are affected with the post-change mean $\mu_{k} = 0.2236$ for all $k,$ we consider two choices of the lower bound $\omega_{k}$'s: one is $\omega_{k} = 0$ for all $k$ and the other is $ \omega_{k} = 0.01$ for all $k.$ The former choice of $\omega_{k} = 0$'s allows us to see the performance of the original SRRS scheme $N_B^{orig}$. For each of these two choices of $\omega_{k}$'s, we vary the linear shrinkage factor $a$ from $0.01$ to $1,$ and then plot the detection delay of $N_{B}$ as a function of $a$ in Figure \ref{figure:Itotbb}. It is clear from Figure \ref{figure:Itotbb} that the linear shrinkage can reduce detection delay from $104.8$ at $a=1$ (baseline $N_B^{orig}$) to $14.0$ at $a=0.17$ when the lower bound $\omega_{k} \equiv 0$ for all $k,$ and can reduce detection delay from $56.6$ at $a=1$  to $11.1$ at $a=0.22$ when the lower bound $\omega_{k} \equiv 0.01$ for all $k.$ Thus both hard-threshold $\omega_{k}$'s and linear shrinkage factor $a$ can reduce the detection delay in this case, though the linear shrinkage factor $a$ seems to be able to play more significant role. This is consistent with our asymptotic results in Section IV.C.

It is interesting to compare the sparse post-change case in  Figure \ref{figure:Itot} with the simultaneous local changes case in Figure \ref{figure:Itotbb}. In both cases, the overall Kullback-Leibler divergence $I_{tot}= 2.5$ are the same, and thus it is not surprising that the original SRRS scheme $N_B^{orig}$ of Lorden and Pollak \cite{lorden2005nonanticipating} has similar detection delays in these two cases (i.e., $104.9$ versus $104.8$). However, the smallest detection delay (i.e., $83.8$) in  Figure \ref{figure:Itot} in the sparse post-change case  is much larger than the smallest detection delay (i.e., $11.1$) in  Figure \ref{figure:Itotbb} when all data streams are affected. In other words, given the same amount of  Kullback-Leibler divergence information, it is much easier to detect simultaneous ``small" local changes in all data streams than to detect ``big" changes in a few unknown data streams if we incorporate relevant prior knowledge appropriately.  This is consistent with our intuition since the latter has to deal with the uncertainty of the subset of affected data streams, which can be very challenging when the dimension $K$ is large.

\subsection{More Simulation About "Curse of Dimensionality"}

In this section, we conduct Monte Carlo simulations to compare the empirical pre-change distributions of the global monitoring statistics $R_n$ in (\ref{new_LR}) under  two different dimensions: $K=1$ and $K=100$, thereby illustrating the challenge of Monte Carlo simulation of the ARL to false alarm when monitoring large $K >1$ number of data streams.

We again assume to monitor $K$ independent normal data streams, and each data stream follows distribution $N(0,1).$ We focus on the performance of the original SRRS scheme $N_B^{orig}$ of Lorden and Pollak \cite{lorden2005nonanticipating}  and the corresponding $R_n$  in (\ref{new_LR}) under the pre-change hypothesis (i.e. $\nu=\infty$). For each scenario of $K=1$ and $K=100,$ and for each time step $n=1,...,1000,$  we ran Monte Carlo to simulate $R_n$ with 2500 replications.

\begin{figure}
  \centering
  \includegraphics[width=2.5in]{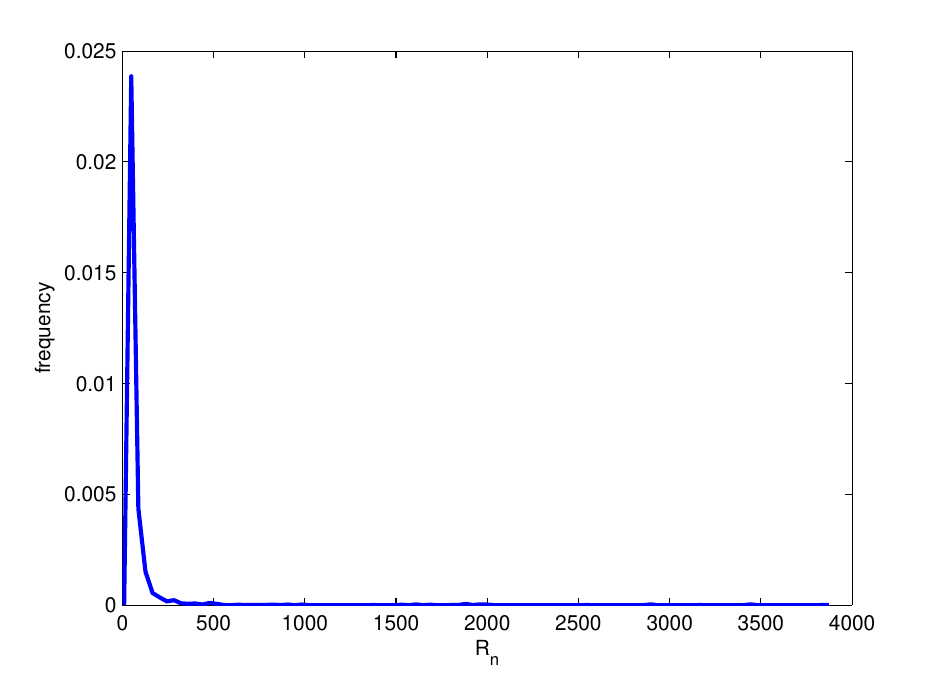} \\
  \includegraphics[width=2.5in]{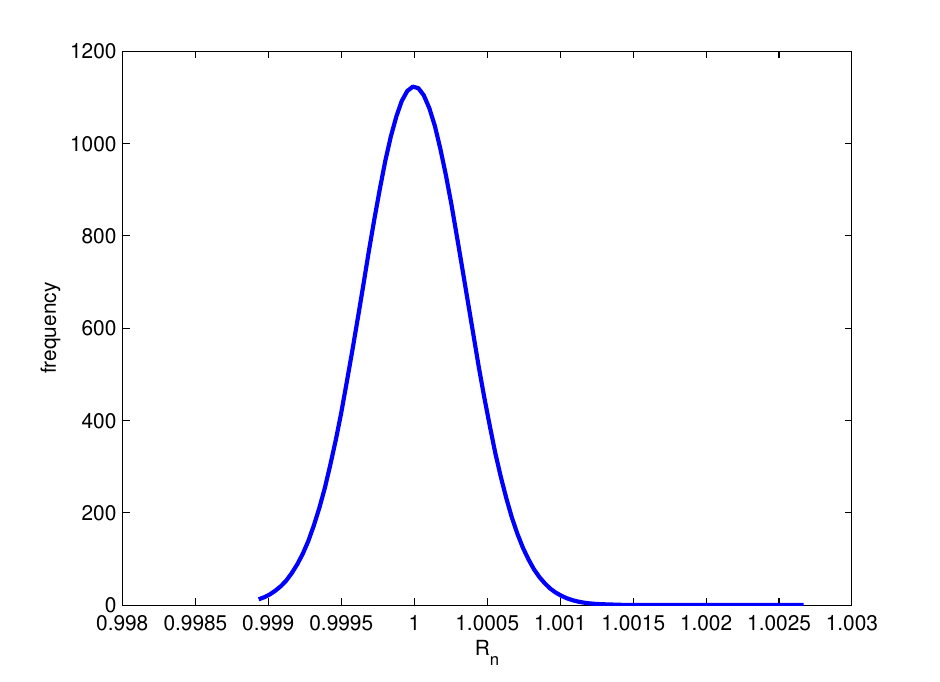}\\
  \caption{Histogram of 2500 simulated $R_n$ with n=500 under two scenarios. {\it Upper Panel:} $K=1,$ and {\it Lower Panel:} $K=100.$} \label{figure:p1hist}
\end{figure}

Figure \ref{figure:p1hist} shows the histogram of $R_n$ at a fixed time $n=500$ for both $K=1$ and $K=100$ cases based on $2500$ replications. As we can see, the empirical distribution of $R_n$ is highly skewed for $K=1$ with values in the range of $[0,4000],$ but $R_{n}$ seems to be empirically normally distributed for $K=100$ with values in the range of $[0.998,1.003].$ Theoretically the empirical mean of $R_{n}$ with $n=500$ should be $500$ no matter whether the dimension $K=1$ or $K=100.$ This suggests that $2500$ replications might be sufficient for $K=1$ dimension, but definitely not large enough for $K=100.$

\begin{figure}
  \centering
  \includegraphics[width=2.5in]{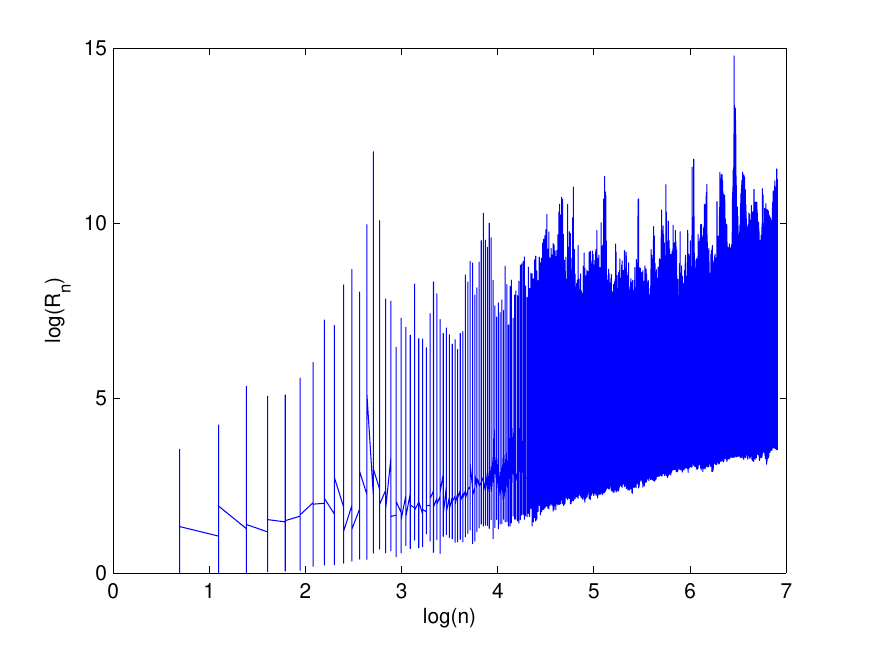} \\
  \includegraphics[width=2.5in]{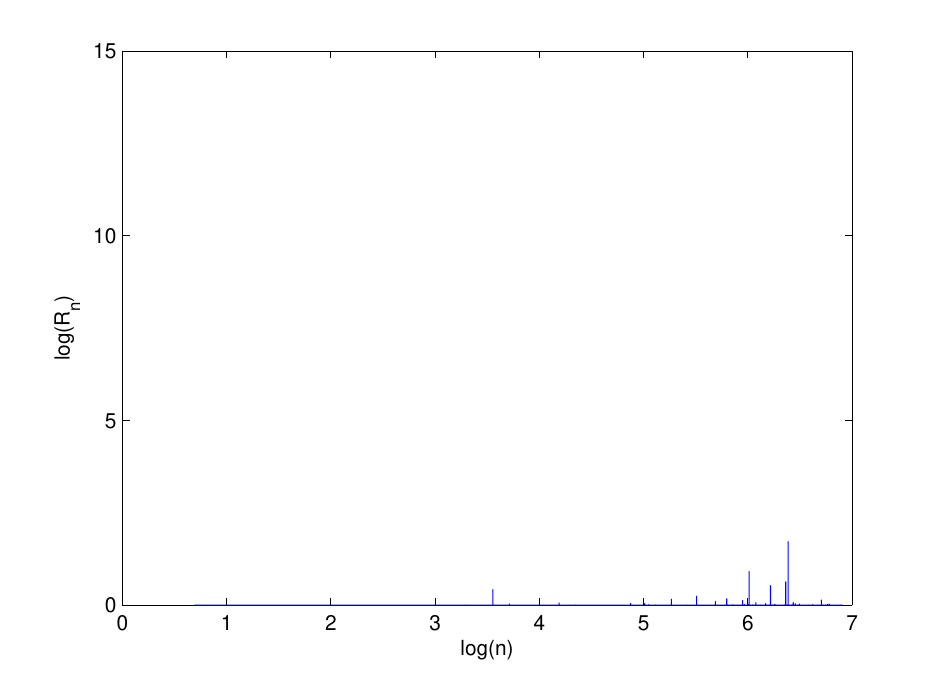}\\
  \caption{
  Plot of $(\log n, \log R_n)$ for $n=1,..,1000$ with 2500 replications under two scenarios. {\it Upper Panel:} $K=1,$ {\it Lower Panel:} $K=100.$} \label{figure: allrep}
\end{figure}

To further explain this issue, we also investigate the dynamic evolution of $R_n$ over time $n.$ To better illustrate, Figure \ref{figure: allrep} plots 2500 simulated $\log(R_n)$ versus $\log(n)$ for both $K=1$ and $K=100$ cases. From Figure \ref{figure: allrep}, there is a clear linear trend of $\log(R_n)$ versus $\log(n)$ when the dimension $K=1,$
which matches the martingale property $\E_{\infty}(R_n) = n.$
On the other hand, when the dimension $K=100,$ most of $\log(R_n)$ are 0, which implies that $2500$ replications are not large enough to represent the property of $R_n.$ The situation is similar even if we increase the number of Monte Carlo runs from $2500$ to a larger number such as $10^4$.
All these simulations results are consistent with our theoretical results.

\section{Conclusions}

In this article, we investigated the quickest change detection problem in the context of monitoring  independent large-scale normally distributed data streams when the post-change means are unknown. The key assumption we make is that  for each individual local data stream, either there are no local changes, or there is a ``big" local change.
Our main contribution is to introduce the shrinkage estimators to quickest change detection, and show that  the shrinkage estimators of the unknown post-change parameters can reduce the overall detection delays by balancing the tradeoff between the first-order and second-order terms of the asymptotic expression on the detection delays.
Specifically, hard-thresholding is attractive in the sparse post-change case when the unknown number of affected data streams is much smaller than the total number of data streams, whereas the linear shrinkage can be useful when all local data streams are affected simultaneously though not necessarily identically. Moreover, we illustrate the challenge of Monte Carlo simulation of
of the ARL to false alarm  when monitoring a large $K$ number of data streams.

While the classical quickest change detection problems have been studied for several decades, further research on the quickest change detection for  monitoring large-scale data streams is needed. For instance, in this paper we focus on the Gaussian model with known variances, and it will be interesting to extend the shrinkage estimators to a more general Gaussian model when the variances in the different data streams are different and unknown under the post-change hypothesis, or to  other distributions such as Poisson. The corresponding theoretical analysis will likely be more challenging, e.g., the definition of $r$ in (\ref{par:par3}) will need to be modified for other distributions such as Poisson. Moreover, it remains an open problem how to overcome the curse of dimensionality to conduct Monte Carlo simulations of the ARL to false alarm efficiently in the context of large-scale data streams. Hopefully this article can stimulate further research on quickest change detection problems in high-dimensional data streams.

\section*{Appendix: Proof of Theorem 2.}

To prove Theorem \ref{thm: detection change}, the crucial technical tools are from those of Theorem 3 and part (iii) of Theorem 4 in  Lorden and Pollak \cite{lorden2005nonanticipating}, which deals with $K=1$ data stream without shrinkage. Below we will highlight the main difference with the dimension $K \ge 1$ and shrinkage.

To simplify the arguments, let us consider the hypothesis testing version of the quickest change detection problem, and assume that we want to test the null hypothesis $H_0:$ {\it no change} against the alternative hypothesis $H_1:$ {\it a change occurs exactly at time $\nu = 1.$} In such a problem,  the corresponding sequential hypothesis testing version of the proposed scheme $N_{B}$ in (\ref{new_LR}) and (\ref{eqnSec03e07})  is defined by
\begin{eqnarray}  \label{eq: RandS}
\tau_{B} =  \inf \{n \geq 1: \Lambda_n \geq B\},
\end{eqnarray}
where the likelihood ratio
\begin{eqnarray} \label{eq: **}
    \Lambda_n = \prod_{\ell=1}^{n}\prod_{k=1}^{K}\frac{f_{\hat{\mu}_{k,\ell}}(X_{k,\ell})}{f_{0}(X_{k,\ell})},
\end{eqnarray}
and the estimate $\hat{\mu}_{k,\ell}$ is a short-handed notation for $\hat{\mu}_{k,m=1,\ell}$ in (\ref{eq: thresh def}).

It is useful to mention that the quickest change detection scheme $N_{B}$ in (\ref{new_LR}) and (\ref{eqnSec03e07}) is closely related to the sequential hypothesis testing procedure $\tau_{B}$ in  (\ref{eq: RandS}) and (\ref{eq: **}), and such a close relation was first discovered in Lorden \cite{lorden1971procedures}. More specifically, for $t=1,2,\ldots,$ denote by $\tau_{B}^{(t)}$ the new stopping time that applies the sequential hypothesis testing procedure $\tau_{B}$ to the data starting from time $t,$ i.e., $\{(X_{1, i}, \ldots, X_{K, i})\}$ for $i=t, t+1, \ldots.$ Then the quickest change detection scheme $N_{B} = \min_{t \ge 1}\{\tau_{B}^{(t)} +t-1\}.$ This relation allows one to show that the detection delay $D(N_{B})$ is asymptotically equivalent to $\E_{\mu}(\tau_{B})$ under the alternative hypothesis $H_1$ when  $\mu=(\mu_1, \cdots, \mu_{K})^{T}$ is the  true post-change mean vector. To emphasize the dependence on the true $\mu,$ denote by $\Prob_{\mu}$ and $\E_{\mu}$ the corresponding probability mean and expectation under the alternative hypothesis $H_1.$ Then it suffices to show that $\E_{\mu}(\tau_{B})$ satisfies the right-hand side of  (\ref{eqn05sec4}).

Recall that in Section IV.B, we denote by $\mu_{k}^{*}$ the limit of $\hat{\mu}_{k,\ell}$ under $\Prob_{\mu}$ for each $k=1, \ldots, K$ as $\ell \rightarrow \infty,$ and define $\mu^* = (\mu_{1}^{*}, \ldots, \mu_{K}^{*})^T.$ A key step of the proof is to relate $\Lambda_n$ in
(\ref{eq: **}) to the likelihood ratio $\Lambda_n^*$ which mis-specify the true post change parameter $\mu_{k}$ of the $k$-th data stream as $\mu_{k}^{*}$ for all $k=1, \ldots, K.$ Since $\log \Lambda_n = 0$ when $n=1,$ we can define the mis-specified log-likelihood ratio by
$$\log \Lambda_n^{*} = \sum_{\ell=2}^{n} \sum_{k=1}^{K} \log \frac{f_{\mu_{k}^{*}}(X_{k,\ell})}{f_{0}(X_{k,\ell})}.$$
for $n \ge 2$ and $\log \Lambda_1^* = 0.$ Then under $\Prob_{\mu},$  $\log \Lambda_n^*$ is a random walk with iid increments that have finite variance and mean $I(\mu^{*}, \mu_0; \mu)$ in  (\ref{eqn06sec3-info}).

For the stopping time $N = \tau_{B}$ in (\ref{eq: RandS}), applying Wald's equation to the random walk $\log \Lambda_n^*$ yields
\begin{eqnarray*}
I(\mu^{*}, \mu_0; \mu) \E_{\mu}(N) &=& \E_{\mu}(\log \Lambda_N^{*}) \cr
 &=& \E_{\mu}(\log \Lambda_N) + \E_{\mu}(\log \Lambda_N^{*} - \log \Lambda_N).
\end{eqnarray*}
For the notational convenience, let $b = \log B.$ Then the standard renewal theorem for over-shoot analysis shows that $\E_{\mu} (\log \Lambda_{N}) = b + O(1)$ for $N = \tau_{B}$ in (\ref{eq: RandS}), where the $O(1)$ term is the over-shoot effect and may depend on the dimension $K$, see Theorem 3 of Lorden and Pollak \cite{lorden2005nonanticipating}.  Thus
\begin{eqnarray}  \label{eqn_0024_wald}
I(\mu^{*}, \mu_0; \mu) \E_{\mu}(N) = b + O(1) + \E_{\mu}(\log \Lambda_N^{*} - \log \Lambda_N).
\end{eqnarray}
Hence it suffices to investigate the property of
\[
\log \Lambda_N^{*} - \log \Lambda_N =  \sum_{\ell=2}^{N} \sum_{k=1}^{K} \log \frac{f_{\mu_{k}^{*}}(X_{k,\ell})}{f_{\hat{\mu}_{k,\ell}}(X_{k,\ell})}
\]
when $N = \tau_{B}$ in (\ref{eq: RandS}) and $X_{k,\ell} \sim N(\mu_{k},1)$ under $\Prob_{\mu}.$

To do so, note that this involves  the likelihood ratio of the form $f_{\mu_{k}^{*}}(X_{k,\ell}) / f_{\phi_{k}}(X_{k,\ell})$ when $X_{k,\ell}$'s are iid $N(\mu_{k},1)$ for each $k,$ and the $\phi_{k}$'s may vary and converge to $\mu_{k}^{*}.$ Thus for any given  $\phi=(\phi_1, \ldots, \phi_{K}),$  we need to define another information number:
\begin{eqnarray} \label{eqn04sec7}
I(\mu^{*}, \phi; \mu) &=&  \E_{\mu} \sum_{k=1}^{K} \left ( \log \frac{f_{\mu_{k}^{*}}(X_{k,\ell})}{f_{\phi_{k}}(X_{k,\ell})} \right )  \cr
&=& \sum_{k=1}^{K} \Big( (\mu_{k}^{*} - \phi_{k}) \mu_{k} - \frac12 (\mu_{k}^{*})^2 +\frac12 (\phi_{k})^2 \Big) \cr
&=& \sum_{k=1}^{K} \Big( (\mu_{k}^{*} - \mu_{k}) \Delta_{k}+ \frac12 (\Delta_{k})^2\Big)
\end{eqnarray}
where $\Delta_{k} = \phi_{k} - \mu_{k}^{*}$ for $k=1, \ldots, K.$ It is useful to compare this new information number with $I(\mu^{*}, \mu_0; \mu)$ in (\ref{eqn06sec3-info1}). On the one hand, they are defined similarly except that $\phi_{k} \equiv \mu_0 =0$ for all $k.$ On the other hand,  $I(\mu^{*}, \mu_0; \mu)$ in (\ref{eqn06sec3-info1}) is related to the first-order term of the detection delay of $\tau_{B},$ whereas $I(\mu^{*}, \phi; \mu)$ in (\ref{eqn04sec7}) contributes to the second-order term of the detection delay when we let $\Delta_{k} = \phi_{k} - \mu_{k}^{*}$ go to $0$ for all $k.$

For any given $\ell= 2, 3, \ldots,$ let $\hat \mu_{\ell}=(\hat \mu_{1,\ell}, \ldots, \hat \mu_{K,\ell})^{T},$ and let $I(\mu^{*}, \hat \mu_{\ell}; \mu)$ be the information number defined in (\ref{eqn04sec7}) when $\phi= \hat \mu_{\ell}$. As in Lorden and Pollak \cite{lorden2005nonanticipating}, the application of the martingale optional sampling theorem to
$\log \Lambda_n^{*} - \log \Lambda_n - \sum_{\ell=2}^{n} I(\mu^{*}, \hat \mu_{\ell}; \mu)$ yields that
\begin{eqnarray} \label{eqn04psec7}
\E_{\mu}\Big(\log \Lambda_N^{*} - \log \Lambda_N \Big) &=& \E_{\mu} \sum_{\ell=2}^{N} I(\mu^{*}, \hat \mu_{\ell}; \mu).
\end{eqnarray}
By (\ref{eqn04sec7}), if we suppress the notation $\ell$ for the sake of convenience and let $\Delta_{k}= \hat{\mu}_{k,\ell} - \mu_{k}^{*},$ then
\begin{eqnarray} \label{eqn027psec7}
\E_{\mu}\Big(I(\mu^{*}, \hat{\mu}_{\ell}; \mu)\Big) &=& \sum_{k=1}^{K} (\mu_{k}^{*} - \mu_{k}) \E_{\mu}(\Delta_{k}) + \cr
&& + \frac12 \sum_{k=1}^{K} \E_{\mu}(\Delta_{k}^2),
\end{eqnarray}
and thus the proof of Theorem \ref{thm: detection change} relies on the  analysis of $\E_{\mu}(\Delta_{k})$ and $\E_{\mu}(\Delta_{k}^2).$

In a high-level description, we may expect that $\Delta_{k}= \hat{\mu}_{k,\ell} - \mu_{k}^{*}$ converges to $0$  as $\ell \rightarrow \infty.$ Hence, for large $\ell,$ we should expect that $\E_{\mu}(\Delta_{k}) \approx 0$ becomes negligible, and the term $\E_{\mu}(\Delta_{k}^2) \approx Var(\Delta_{k})$ may or may not be significant. Indeed,  for a given $\ell,$ we will show below that as $\ell \rightarrow \infty,$
\begin{eqnarray} \label{eqn034psec7}
\E_{\mu}(\Delta_{k}) = o(\frac{1}{(\ell-1)^2})
\end{eqnarray}
and
\begin{eqnarray} \label{eqn035psec7}
\E_{\mu}(\Delta_{k}^2) \sim
\left\{
  \begin{array}{ll}
    a^2/ (\ell-1), & \hbox{ if $|\mu_{k}| > \omega_{k}$;} \\
    o(\frac{1}{(\ell-1)^2}), & \hbox{ if $|\mu_{k}| < \omega_{k}.$}
  \end{array}
\right.
\end{eqnarray}
Let us postpone the proof of (\ref{eqn034psec7}) and (\ref{eqn035psec7}) in a little bit, and apply them directly  to (\ref{eqn027psec7}),  we have
\begin{eqnarray*}
& & \E_{\mu}\Big(I(\mu^{*}, \hat{\mu}_{\ell}; \mu)\Big)\\
 &=& \sum_{k=1}^{K} o\Big(\frac{1}{(\ell - 1)^2} \Big)
+ \frac12  \sum_{k: |\mu_{k}| > \omega_{k}} \frac{a^2}{\ell-1}  + \cr
&& + \frac12 \sum_{k: |\mu_{k}| < \omega_{k}} o\big(\frac{1}{(\ell-1)^2} \big) \Big) \\
&=&  \frac{r}{2} \frac{a^2}{\ell - 1} + o(\frac{1}{(\ell - 1)^2})
\end{eqnarray*}
as $\ell$ goes to $\infty,$ where $r$ is defined in (\ref{par:par3}).  Plugging this into  (\ref{eqn04psec7}), we have
\[
\E_{\mu}\Big(\log \Lambda_N^{*} - \log \Lambda_N \Big) 
= \frac{r a^2}{2} (1+o(1)) \E_{\mu} \sum_{\ell=2}^{N} \frac{1}{(\ell - 1)}.
\]
The summation of the above relation can then be estimated as in  Theorem 3 of Lorden and Pollak \cite{lorden2005nonanticipating} by
\[
(1+o(1)) \sum_{\ell=2}^{n_0} \frac{1}{(\ell - 1)} \approx (1+o(1)) \log(n_0)
\]
where $n_0 =$ the largest integer $\le \E_{\mu}(N).$ Combining this with (\ref{eqn_0024_wald}) yields
\begin{eqnarray*}
& & I(\mu^{*}, \mu_0; \mu) \E_{\mu}(N) \\
&=& b + O(1) + \E_{\mu}(\log \Lambda_N^{*} - \log \Lambda_N)  \\
&=& b + O(1) + (1+o(1)) \frac{r a^2}{2} \log(\E_{\mu}(N)).
\end{eqnarray*}
This gives an equation for $\E_{\mu}(N),$ and thus $\E_{\mu}(N)$ can be found by solving the equation of the form $x = \alpha + \beta \log(x)$ for large $\alpha > 0$ and possibly large $\beta > 0.$
Taking logarithms of both sides yields
\begin{eqnarray*}
\log(x) &=& \log(\alpha + \beta \log(x)) = \log \max\{\alpha, \beta \log(x)\} + O(1) \cr
&=& \max\{\log \alpha, \log \beta\} + o(\log x),
\end{eqnarray*}
where we use the fact that $\max(x, y) \le x+y \le 2 \max(x,y)$ for $x > 0, y > 0$ and $O(1) = O(\log \log x) = o(\log x)$ for large $x.$ Plugging this relation back to $x = \alpha + \beta \log(x)$
yields that $$x = \alpha + (1+o(1)) \beta  \max\{\log \alpha, \log \beta\}.$$
Using the above arguments to derive $\E_{\mu}(N)$ and absorbing all insignificant terms to the $o(1)$ term,  we have
\begin{eqnarray*}
& & \E_{\mu} (N) \\
&=& \left ( b +  (1+o(1)) \frac{r a^2}{2}  \log\frac{\max\{b, r a^2/2\}}{I(\mu^{*}, \mu_0; \mu)} \right )/ I(\mu_*,\mu_0; \mu )\\
\end{eqnarray*}
which becomes the right-hand side of (\ref{eqn05sec4}) as $b= \log(B)$ goes to $\infty.$  Thus the theorem holds.

It remains to prove (\ref{eqn034psec7}) and (\ref{eqn035psec7}). The details can be simplified to the following elementary probability question. Given two real numbers $\mu$ and $\omega > 0,$ and $|\mu| \ne \omega.$ Assume $Y = (X_{k,1} + \ldots + X_{k, \ell-1}) / (\ell - 1) \sim N(\mu, \sigma^2 = 1 /(\ell - 1)),$ and define a new random variable
\[
Y^{*} = \left\{
          \begin{array}{ll}
            aY + b, & \hbox{if $|Y| \ge \omega$;} \\
            c, & \hbox{if $|Y| < \omega.$}
          \end{array}
        \right.
\]
and a new constant
\[
\mu^{*} = \left\{
          \begin{array}{ll}
            a \mu + b, & \hbox{if $|\mu| > \omega$;} \\
            c, & \hbox{if $|\mu| < \omega.$}
          \end{array}
        \right.
\]
Let $\Delta = Y^{*} - \mu^{*},$ and  we want to show the asymptotic properties of $\E(\Delta)$  and $\E(\Delta^2)$  satisfy (\ref{eqn034psec7}) and (\ref{eqn035psec7}) as $\sigma^2 = \frac{1}{\ell - 1} \rightarrow 0.$

We need to consider three cases, depending on the relationship between $\mu$ and $\pm \omega.$ Following the traditional notation, let $Z = (Y - \mu) / \sigma \sim N(0,1),$ and denote by $\phi(z)$ and $\Phi(z)$ for the probability density function (pdf) and cumulative distribution function (cdf) of $N(0,1),$ respectively. Also define $\lambda_1 = (-\omega - \mu) / \sigma$ and $\lambda_2 = (\omega - \mu) / \sigma.$ Then $Y^{*} = (a (\mu + \sigma Z) + b) (1\{ Z \le \lambda_1\}  + 1\{ Z \ge \lambda_2\}) + c  (1\{\lambda_1 \le Z \le \lambda_2\}).$

Let us focus on the case when $\mu > \omega.$ In this case, we have $\mu^{*} = a \mu + b$ and $\lambda_1 < \lambda_2 \rightarrow -\infty$ as $\sigma \rightarrow 0.$  Hence,
\begin{eqnarray*}
\Delta = Y^{*} - \mu^{*} &=& a \sigma Z (1\{ Z \le \lambda_1\}  + 1\{ Z \ge \lambda_2\}) \cr
&&  + (c - a \mu - b) 1\{\lambda_1 < Z < \lambda_2\}.
\end{eqnarray*}
Since $\lambda_1 < \lambda_2 \rightarrow -\infty$ as $\sigma \rightarrow 0,$  the event $1\{ Z \ge \lambda_2\}$ is dominant whereas the other two events are rare events. Thus we should expect that $\Delta \approx a \sigma Z,$ and thus $\E(\Delta) \approx o(\sigma^4)$ and $\E(\Delta^2) \approx Var(a \sigma Z) = a^2 \sigma^2.$ To be more rigorous,
\begin{eqnarray*}
\E(\Delta) &=& \int_{-\infty}^{\lambda_1} a \sigma z \phi(z) d z + \\
&& +  \int_{\lambda_2}^{\infty} a \sigma z \phi(z) d z + \int_{\lambda_1}^{\lambda_2} (c -a\mu-b) \phi(z) d z \\
&=& - a \sigma \phi(|\lambda_1|) +  a \sigma \phi(|\lambda_2|) + \\
&&  + (c -a\mu-b) \Prob(|\lambda_2| \le Z \le |\lambda_1| )
\end{eqnarray*}
where we use the fact $\int_{-\infty}^{\lambda} z \phi(z) d z = - \phi(|\lambda|) = - \int_{\lambda}^{\infty} z \phi(z) d z$ when $\lambda < 0.$
By the well-known fact that $\frac{x}{1 + x^2} \phi(x) \le \Prob( Z > x) \le \frac{\phi(x)}{x}$ for all $x \ge 0,$ it is clear that $\E(\Delta) = O( \sigma \phi(|\lambda_1|)) + O(\sigma \phi(|\lambda_2|)) = o(\sigma^4)$ as $\sigma$ goes to $0,$ since $O(\phi(x/\sigma)) = O(\exp(- \frac{x^2}{2 \sigma^2})) = o(\sigma^4)$ for any $x \ne 0.$

In addition,
\begin{eqnarray*}
\E(\Delta^2) &=& \int_{-\infty}^{\lambda_1} (a \sigma z)^2 \phi(z) d z + \\
  && + \int_{\lambda_2}^{\infty} (a \sigma z)^2 \phi(z) d z + \int_{\lambda_1}^{\lambda_2} (c -a\mu-b)^2 \phi(z) d z \\
&=& a^2 \sigma^2 [\Prob(Z > |\lambda_1|) + |\lambda_1|\phi(|\lambda_1|)] \\
 && + a^2 \sigma^2 [1- \Prob(Z > |\lambda_2|) + |\lambda_2|\phi(|\lambda_2|)] \\
 && + (c -a\mu-b)^2 \Prob(|\lambda_2| \le Z \le |\lambda_1| )\\
&=& a^2 \sigma^2 + o(\sigma^4).
\end{eqnarray*}
Here in the second equation, we use the fact that $\int_{-\infty}^{\lambda} z^2 \phi(z) d z = \Prob(Z > |\lambda|) + |\lambda|\phi(|\lambda|) = 1- \int_{\lambda}^{\infty} z^2 \phi(z) d z $ for $\lambda < 0,$ which follows from the integration by parts for $z^2 \phi(z) = - z (\phi(z))^{'}.$  Thus (\ref{eqn034psec7}) and (\ref{eqn035psec7}) hold for the case when  $\mu > \omega.$

The above arguments can be easily extend to the other cases when $\mu < - \omega$ or $- \omega < \mu < \omega.$ For instance, when $- \omega < \mu < \omega,$ we have $\mu^* = c$ and $\lambda_1 \rightarrow -\infty, \lambda_2 \rightarrow \infty$ as $\sigma \rightarrow 0.$ Thus $\Delta = Y^{*} - c = (a (\mu + \sigma Z) + b - c) (1\{ Z \le \lambda_1) + 1\{Z \ge \lambda_2\}).$  Since the probabilities of both events $1\{Z \le -\lambda_1\}$ and $1\{Z \ge \lambda_2\}$  go to $0$ exponentially as $\sigma$ goes to $0,$ the above arguments can show that both $\E(\Delta)$ and $\E(\Delta^2)$ are negligible (order $o(\sigma^4)$), completing the proof of the theorem.


\section*{Acknowledgement}
This work was supported in part by the NSF grants DMS-0954704 and CMMI-1362876.
The authors are grateful to the Associate Editor and two anonymous reviewers for the detailed and constructive comments that greatly improve the quality and presentation of the article.
\par



\end{document}